\documentclass [11pt,oneside]{article}

\reversemarginpar \textwidth=14cm \textheight=19cm

\usepackage{amssymb} \usepackage{amsfonts} \usepackage{amsmath}
\usepackage{amsthm} \usepackage{epsfig}
\usepackage{wrapfig}

\newcommand{\mettifig}[1]{\epsfig{file=#1}}

\newcommand{\altramettifig}[1]{\includegraphics[width=\textwidth]{#1}}

\newcommand{\lens}[2]{L({\scriptstyle #1},{\scriptstyle #2})}
\newcommand{\seifeul}[2]{\big(#1,#2\big)}
\newcommand{\seifuno}[3]{\big(#1,({\scriptstyle #2},{\scriptstyle #3})\big)}
\newcommand{\seifunoeul}[4]{\big(#1,({\scriptstyle #2},{\scriptstyle #3}),#4\big)}
\newcommand{\seifdue}[5]{\big(#1,({\scriptstyle #2},{\scriptstyle #3}),
                       ({\scriptstyle #4},{\scriptstyle #5})\big)}
\newcommand{\seifdueeul}[6]{\big(#1,({\scriptstyle #2},{\scriptstyle #3}),
                       ({\scriptstyle #4},{\scriptstyle #5}),#6\big)}
\newcommand{\seiftre}[7]{\big(#1,({\scriptstyle #2},{\scriptstyle #3}),
                       ({\scriptstyle #4},{\scriptstyle #5}),
                       ({\scriptstyle #6},{\scriptstyle #7})\big)}
\newcommand{\seiftreeul}[8]{\big(#1,({\scriptstyle #2},{\scriptstyle #3}),
                       ({\scriptstyle #4},{\scriptstyle #5}),
                       ({\scriptstyle #6},{\scriptstyle #7}),#8\big)}
\newcommand{\bigu}[4]{\bigcup\nolimits_{{\tiny{\matr {#1} {#2} {#3} {#4}}}\phantom{\Big|}\!\!}}
\newcommand{\bigb}[4]{\big/_{{\tiny{\matr {#1} {#2} {#3} {#4}}}\phantom{\Big|}\!\!}}

\newtheorem{mainteo}{Theorem}
\newtheorem{lemma}{Lemma}[section] 
\newtheorem{teo}[lemma]{Theorem}
\newtheorem{rem}[lemma]{Remark} 
\newtheorem{prop}[lemma]{Proposition}
\newtheorem{cor}[lemma]{Corollary}

\newcommand{\matN}{\ensuremath {\mathbb{N}}}
\newcommand{\matR} {\ensuremath {\mathbb{R}}}
\newcommand{\matQ} {\ensuremath {\mathbb{Q}}}
\newcommand{\matZ} {\ensuremath {\mathbb{Z}}}
\newcommand{\matC} {\ensuremath {\mathbb{C}}}

\newcommand{\matH} {\ensuremath {\mathbb{H}}}
\newcommand{\matRP} {\ensuremath {\mathbb{RP}}}

\newcommand{\calG}{\ensuremath {\mathcal{G}}}

\newcommand{\nota} [1] {\caption{\footnotesize{#1}}}
\newcommand{\matr} [4] {\left(\begin{array}{@{}c@{\ }c@{}} #1 & #2 \\ #3 & #4 \\ \end{array} \right)}

\newfont{\Got}{eufm10 scaled 1200}

\newcommand{\GL}{{\rm GL}}

\font\titsc=cmcsc10 scaled 1200

\newcommand{\dimo}[1]{\vspace{2pt}\noindent\textit{Proof of \ref{#1}}.\ }
\newcommand{\finedimo}{{\hfill\hbox{$\square$}\vspace{2pt}}}

\newcommand{\timtil}{\begin{picture}(12,12)
\put(2,0){$\times$}\put(2,4.5){$\sim$}\end{picture}}

\newcommand{\grafetto}{
\begin{picture}(21,8)
\put(4.2,4){\circle{8}}
\put(8,4){\circle*{2}}
\put(8,4){\line(1,0){8}}
\put(16,4){\circle*{2}}
\end{picture}
}

\author{Bruno \titsc{Martelli}
\and Carlo \titsc{Petronio}}

\title{Dehn filling of the ``magic'' 3-manifold}

\begin{document}

\maketitle

\begin{abstract}
  We classify all the non-hyperbolic Dehn fillings of the complement of the
  chain-link with 3 components, conjectured to be the smallest
  hyperbolic 3-manifold with 3 cusps. We deduce the classification of
  all non-hyperbolic Dehn fillings of infinitely many 1-cusped and 2-cusped
  hyperbolic manifolds, including most of those with smallest known volume.

Among other consequences of this classification, we mention the
following:

\begin{itemize}
\item For every integer $n$ we can prove that there are infinitely
many hyperbolic knots in $S^3$ having exceptional surgeries $\{n,
n+1, n+2, n+3\}$, with $n+1, n+2$ giving small Seifert manifolds
and $n, n+3$ giving toroidal manifolds; \item We exhibit a
2-cusped hyperbolic manifold that contains a pair of inequivalent
knots having homeomorphic complements; \item We exhibit a chiral
$3$-manifold containing a pair of inequivalent hyperbolic knots
with orientation-preservingly homeomorphic complements; \item We
give explicit lower bounds for the maximal distance between small
Seifert fillings and any other kind of exceptional filling.
\end{itemize}

  \vspace{4pt}

\noindent MSC (2000): 57M27 (primary), 57M20, 57M50 (secondary).
\end{abstract}

\begin{wrapfigure}{r}{1.7in}
  \centering
  \mettifig{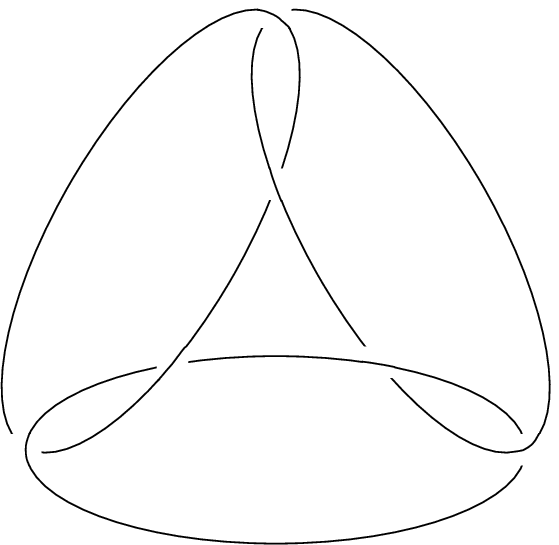, width = 3 cm}
  \nota{The chain link.}
  \label{chain:fig}
\end{wrapfigure}

We study in this paper the Dehn fillings of the complement $N$ of
the chain-link with 3 components in $S^3$, shown in
Fig.~\ref{chain:fig}.

The hyperbolic structure of $N$ was first constructed by Thurston
in his notes~\cite{bibbia}, and it was also noted there that the
volume of $N$ is particularly small. The relevance of $N$ to
3-dimensional topology comes from the fact that filling $N$ one
gets most of the hyperbolic manifolds known and most of the
interesting non-hyperbolic fillings of cusped hyperbolic
manifolds. For these reasons $N$ was called the ``magic manifold''
by Gordon and Wu~\cite{GoWu-TA, gordon}. It appears as $M6_1^3$
in~\cite{Ca-Hi-We} and it is the
hyperbolic manifold with $3$ cusps of smallest known volume and of
smallest complexity~\cite{AdSh}. (We refer here to the complexity
defined by Matveev in~\cite{Mat}, and we mean that $N$ has an
ideal triangulation with $6$ tetrahedra, while all other
hyperbolic manifolds with $3$ cusps need more than $6$
tetrahedra).

The main results of this paper, stated in
Section~\ref{main:section}, provide a complete classification
(with prime and JSJ decomposition) of all non-hyperbolic Dehn
fillings of $N$, including partial ones. Quite some energy has
been devoted in the literature to the understanding of
\emph{exceptional} slopes on hyperbolic manifolds,
\emph{i.e.}~slopes giving non-hyperbolic fillings, but a complete
description of the actual filled manifolds is seldom available.
Our tables provide a complete description of the exceptional
fillings for many hyperbolic manifolds with $1$ and $2$ cusps: for
this reason we hope that our results will prove useful. We now
list some consequences of this classification.

\paragraph{Surgery along knots in $S^3$}
If the complement of a knot arises by filling two cusps of $N$,
its exceptional surgeries can be recovered from our tables. Every
such knot has at least $4$ exceptional surgeries. In particular,
we prove the following:

\begin{mainteo} \label{Seifert:teo}
For every integer $n$ there are infinitely many hyperbolic knots in $S^3$ whose
exceptional surgeries are $\{n,n+1,n+2,n+3\}$, with $n$ and $n+3$ giving
toroidal manifolds, and $n+1, n+2$ giving small
Seifert manifolds.
\end{mainteo}

Motegi and Song recently proved  that for every
$n$ there is a hyperbolic knot whose $n$-surgery is small
Seifert~\cite{MoSo}, and Teragaito
proved the same result for toroidal surgeries~\cite{Te}.
Theorem~\ref{Seifert:teo} exhibits infinitely many examples for
each $n$. Other knots in $S^3$ with Seifert surgeries are also
described in~\cite{De,EuMu,MaMiMo}.

\paragraph{Seifert fillings}
Following Gordon~\cite{gordon}, let us say that a manifold is of
type $D$, $A$, $S$, or $T$ if it contains respectively an
essential disc, annulus, sphere, or torus, and of type $S^H$ or
$T^H$ if it contains a Heegaard sphere or torus.
An exceptional filling of a hyperbolic
$3$-manifold is either of one of these types, or a small closed
Seifert manifold, or a counterexample to the Geometrization
Conjecture (which seems ruled out by Perelman). Detecting small
Seifert fillings is the most difficult task, since these manifolds
do not contain any small essential embedded surface (they only
contain \emph{immersed} essential tori).

Denoting by $\Delta$ the minimal geometric intersection between slopes on a torus, for
$X_0,X_1\in\{D,A,S,T,S^H,T^H\}$ we define as in~\cite{gordon} the number
$\Delta(X_0,X_1)$ as the maximum
of $\Delta(\alpha_0,\alpha_1)$ over $M,\alpha_0,\alpha_1$, where $M$ is cusped
hyperbolic, $\alpha_0,\alpha_1$ are slopes on the same cusp of $M$, and
$M(\alpha_i)$ is of type $X_i$.
All values of $\Delta(X_0,X_1)$ are known, except $\Delta(S,S^H)$ and $\Delta(T,T^H)$,
which are conjectured to be respectively $-\infty$ (this is the Cabling Conjecture) and $3$.
Our tables support both these conjectures.

Let us now
say that a small Seifert manifold is of type $Z$. As opposed to
most values of $\Delta(X_0,X_1)$ for
$X_0,X_1\in\{D,S,A,T,S^H,T^H\}$, no value $\Delta(X,Z)$ seems to
be known for $X\in\{S,T,S^H,T^H,Z\}$. Our tables imply the
following lower bounds:

\begin{center}
\begin{tabular}{c||c|c|c|c|c}
    & $S$ & $T$ & $S^H$ & $T^H$ & $Z$ \\
\hline\hline $Z$ & $2$ & $7$ & $1$& $2$ & $6$
\end{tabular}
\end{center}

\emph{Infinitely} many hyperbolic manifolds realize the following values:
\begin{center}
\begin{tabular}{c||c|c|c|c|c}
    & $S$ & $T$ & $S^H$ & $T^H$ & $Z$ \\
\hline\hline $Z$ & $2$ & $5$ & $1$   & $2$ & $4$
\end{tabular}
\end{center}

We have not found in the literature examples
realizing greater values than those given by these tables.

\paragraph{Inequivalent knots with homeomorphic complements}
The celebrated result of Gordon and Luecke~\cite{GoLu} states that
inequivalent knots in $S^3$ have non-homeomorphic complements.
Here, two knots are equivalent when there is a self-homeomorphism
of the ambient $3$-manifold sending one to the other. Bleiler,
Hodgson, and Weeks found two inequivalent hyperbolic knots in
$\lens{49}{18}$ with homeomorphic complements (but
the homeomorphism reverses the orientation). This
example can be seen in our tables, and is still the only one known
inside a closed manifold.

More examples arise in manifolds with
boundary. There are infinitely many pairs of inequivalent
hyperbolic knots with homeomorphic complements in the solid torus
(their construction, due to Berge~\cite{Be} and Gabai~\cite{Ga2},
led to the example in $\lens{49}{18}$ mentioned above) and in
handlebodies of arbitrary genus~\cite{FriMaPe3}. We find in this
paper some examples inside other manifolds:

\begin{mainteo} \label{knots1:teo}
The Whitehead sister link complement contains a pair
of inequivalent hyperbolic knots with homeomorphic complements.
\end{mainteo}

\begin{mainteo} \label{knots2:teo}
There exists a chiral Seifert manifold fibred over the disc with two
exceptional fibres that contains a pair of inequivalent hyperbolic
knots with orientation-preservingly homeomorphic complements.
\end{mainteo}

The example provided by Theorem~\ref{knots1:teo} is probably the
first one discovered in a hyperbolic manifold. We note that the
two inequivalent knots are geodesics of distinct length both
having $N$ as a complement (but the homeomorphism between the
complements reverses the orientation). The example provided by
Theorem~\ref{knots2:teo} is probably the first one of knots with
orientation-preservingly homeomorphic complements in a chiral
manifold (all other known examples are inside handlebodies).

\paragraph{Structure of the paper}
The list of all non-hyperbolic Dehn
fillings of $N$ is presented and discussed in
Section~\ref{main:section}, where we state several theorems and
deduce from them some of the results announced in the present
Introduction, and a few more. Section~\ref{proofs:section}
contains the proofs of the classification theorems. The Appendix
contains more tables (deduced from the lists in
Section~\ref{main:section}) showing the non-hyperbolic fillings of
some notable classes of manifolds with one or two cusps. These
further tables are used to prove the rest of the results stated
above.

We warmly thank Jeff Weeks for very helpful suggestions on how to use SnapPea
for the last step of the proof of Theorem~\ref{three:cusps:teo}.

\section{Non-hyperbolic fillings of the chain-link complement} \label{main:section}

We state here the main results of the paper, namely Theorems~\ref{one:cusp:teo},~\ref{two:cusps:teo},
and~\ref{three:cusps:teo} which classify all non-hyperbolic fillings of the chain-link complement $N$.
We then use these results to prove the theorems stated in the Introduction.

\subsection{Main results}

To give our statements, we first make some conventions.

\paragraph{Homology bases}
We fix the standard (meridian, longitude) homology basis on each component of
$\partial N$. Using it, we identify a slope $\alpha$ with a number
in $\matQ\cup\{\infty\}$, and we denote by $N(\alpha)$ the Dehn filling of $N$
along $\alpha$. Any two components of $\partial N$ can be interchanged
by an automorphism which preserves the bases, so $N(\alpha)$
depends on $\alpha$ only. On each partial filling
$N(\alpha)$ or $N(\alpha,\beta)$ we take the induced homology bases,
and we note that the two components of $N(\alpha)$ are still symmetric.

\paragraph{Lens spaces} For simplicity we employ
a more flexible notation for lens spaces than the usual one. We write
$\lens 2q$ for $\matRP^3$, $\lens 1q$ for $S^3$, $\lens 0q$ for
$S^2\times S^1$, and $\lens pq$ for $\lens {|p|}{q'}$ with $q\equiv
q'\ (\textrm{mod } p)$ and $0 < q' < |p|$, for any coprime $p,q$.

\paragraph{Torus bundles} We denote by $T_X$ the torus bundle with
monodromy $X\in\GL_2(\matZ)$.

\paragraph{Seifert manifolds}
If $F$ is a surface, $k\in\matN$, $b\in\matZ$ and $\{(p_i,q_i)\}_{i=1}^k$
are coprime pairs with $|p_i|\geqslant 2$, we define a Seifert manifold
$\big(F, ({\scriptstyle p_1},{\scriptstyle q_1}),\ldots,({\scriptstyle p_k},{\scriptstyle q_k}),b\big)$
as follows. Let $\Sigma$ be $F$ minus $k+1$ open discs and let $W$
be the orientable $S^1$-bundle over $\Sigma$. Give $W$ any orientation,
pick a section $\sigma$ of $W$, and choose positive homology bases
$\{(\mu_i,\lambda_i)\}_{i=0}^k$ on the components of $\partial W$ arising
from the punctures of $F$, with $\mu_i\subset\partial\sigma$
and a fibre as $\lambda_i$. Then
$\big(F, ({\scriptstyle p_1},{\scriptstyle q_1}),\ldots,({\scriptstyle p_k},{\scriptstyle q_k}),b\big)$
is the Dehn filling of $W$ along
$p_1\mu_1+q_1\lambda_1,\ldots,p_k\mu_k+q_k\lambda_k,\mu_0+b\lambda_0$.

We say the parameters are \emph{normalized} if $p_i>q_i>0$. In our statements
we will allow parameters to be non-normalized only when they depend on some
variable. Note that $b$ can be dismissed when $F$ has boundary and that
our ``filling parameters'' $(p_i,q_i)$ are not the usual Seifert orbital parameters
$(a_i,b_i)$, which are computed as $a_i=|p_i|$ and
$b_i\equiv q_1^{-1}\cdot {p_i}/{|p_i|} \ ({\rm mod}\ |p_i|)$.

We denote by $D$, $A$, $S$, $P$, and $T$
respectively the disc, the annulus, the M\"obius strip, the pair of pants, and the torus.

\paragraph{Graph manifolds}
If $F$ is an orientable surface with $h\geqslant 1$ boundary components
and $\Sigma$ is $F$ minus $k$ discs, we proceed as above to construct
homology bases $\{(\mu_i,\lambda_i)\}_{i=1}^{h+k}$ on
$\partial(\Sigma\times S^1)$, with the extra requirement that
the $\mu_i$'s should be oriented as components of $\partial\Sigma$.
For coprime pairs $\{(p_i,q_i)\}_{i=1}^k$ with $|p_i|\geqslant 2$
we then get a Seifert manifold
$\big(F, ({\scriptstyle p_1},{\scriptstyle q_1}),\ldots,
({\scriptstyle p_k},{\scriptstyle q_k})\big)$ with fixed homology bases
on its $h$ boundary components. If $M,M'$ are Seifert manifolds
with boundary and orientable base surfaces, and $X\in\GL_2(\matZ)$, we can
then define without ambiguity
$M\cup_X M'$ to mean $M\cup_f M'$ with $f:T\to T'$, where $T$ and
$T'$ are arbitrary boundary tori of $M$ and $M'$, respectively,
and $f$ acts on homology as $X$ with respect to the fixed bases.
Similarly we can define $M\big/_X$ when $M$ has at least 2 components.

\bigskip

We are now ready to give our statements, which correspond to the
filling of 1, 2, or all 3 components of $\partial N$,
respectively. Throughout this paper `hyperbolic' means
`finite-volume complete hyperbolic,' and the symbol $N$
is only used to denote the complement of the 3-components chain-link in
$S^3$, shown in Fig.~\ref{chain:fig}.

\begin{teo}\label{one:cusp:teo}
$N(\alpha)$ is hyperbolic unless $\alpha\in\{\infty,-3,-2,-1,0\}$, and:
\begin{itemize}
\item $N(\infty) = T\times I$;
\item $N(-3) = \seifuno A21 \bigu 0110 \seifuno A21$;
\item $N(-2) = \seifuno P31$;
\item $N(-1) = \seifuno P21$;
\item $N(0) = \seifdue D2131 \bigu 11{-1}0 \big(P\times S^1\big)$.
\end{itemize}
\end{teo}

\begin{teo}\label{two:cusps:teo}
$M=N(p/q, r/s)$ is hyperbolic except
if one of the following occurs up to permutation:
\begin{itemize}
\item
$p/q\in\{\infty,-3,-2,-1,0\}$;
\item
$(p/q,r/s)\in\{(1,1),(-4,-1/2),(-3/2,-5/2)\}$.
\end{itemize}
In all cases
$M$ is as shown in Table~\ref{two:cusps:table},
where $n$ stands for a generic integer and $r/s\neq
\infty$ unless $p/q=\infty$.
\end{teo}

\begin{table}\begin{center}
\begin{tabular}{c|c||l}
$p/q$ & $r/s$ & \qquad $M$ \\
        \hline\hline
$\infty$ & any & $D\times S^1$ \phantom{\Big|}  \\
        \hline
$-3$
&
    $-1$ &
    $\matRP^3\#(D\times S^1)$ \phantom{\Big|} \\
        \cline{2-3}
&
    $-2$ &
    $D\times S^1$ \phantom{\Big|}\\
        \cline{2-3}
&
    $-3$ &
    $\seifuno S21$ \phantom{\Big|}\\
        \cline{2-3}
&
    $-1+\frac1n, \ne -2$ &
    $\seifdue D21{2n+1}2$ \phantom{\Big|} \\
        \cline{2-3}
&
    $\ne -1, \ne -1+\frac1n$  &
    $\seifdue D21{r+s}s\bigu 0110 \seifuno A21$ \phantom{\Big|} \\
        \hline
$-2$
&
    $-2$ &
    $\lens 31\#(D\times S^1)$ \phantom{\Big|} \\
        \cline{2-3}
&
    $-2+\frac1n$ &
    $D\times S^1$ \phantom{\Big|}\\
        \cline{2-3}
&
    $\ne -2, \ne -2+\frac1n$  &
    $\seifdue D31{r+2s}s$ \phantom{\Big|} \\
        \hline
$-1$
&
    $-3$ &
    $\matRP^3\#(D\times S^1)$ \phantom{\Big|} \\
        \cline{2-3}
&
    $-3+\frac1n$ &
    $D\times S^1$ \phantom{\Big|}\\
        \cline{2-3}
&
    $\ne -3, \ne -3+\frac1n$  &
    $\seifdue D21{r+3s}s$ \phantom{\Big|} \\
        \hline
$0$
&
    $\in\matZ$ &
    $\seifdue D2131$ \phantom{\Big|} \\
        \cline{2-3}
&
    $\not\in\matZ$ &
    $\seifdue D2131 \bigu 11{-1}0 \seifuno Asr$ \phantom{\Big|} \\
        \hline
$1$ & $1$ & $\big(P\times S^1\big)\bigb 0{1}{1}0$ \phantom{\Big|} \\
        \hline
$-\frac32$ & $-\frac52$ & $\seifdue D2131 \bigu 11{-1}0 \seifuno A21$ \phantom{\Big|} \\
        \hline
$-4$ & $-\frac 12$ & $\seifdue D2131 \bigu 11{-1}0 \seifuno A21$ \phantom{\Big|}
\end{tabular}
\nota{Classification of non-hyperbolic fillings of $N$ with 1 boundary component}\label{two:cusps:table}
\end{center}\end{table}

\begin{teo}\label{three:cusps:teo}
$M=N(p/q,r/s,t/u)$ is hyperbolic except
if one the following occurs up to permutation:
\begin{itemize}
\item $p/q=\infty$, in which case $M=\lens{tr-us}{tr'-us'}$, where
$r', s'\in\matZ$ are such that $rs'-sr'=\pm1$;
\item $p/q\in\{-3,-2,-1,0\}$, in which case
$M$ is as described in Table~\ref{closed:fill:table1}
(with $n$ and $m$ standing for generic integers);
\item $(p/q,r/s)\in\{(1,1),(-4,-1/2),(-3/2,-5/2)\}$, in which case
$M$ is as described in Table~\ref{closed:fill:table2};
\item $(p/q,r/s,t/u)$ is one of the 14 triples listed in Table~\ref{closed:fill:table3},
together with a description of the corresponding $M$.
\end{itemize}
\end{teo}

\begin{table}\begin{center}
\begin{tabular}{c||c|c||l}
$p/q$ & $r/s$ & $t/u$ & \qquad $M$ \\
        \hline\hline
$-3$
     &
    $1$ or $-\frac53$ &
    $r/s$ &
    $(K,1)\phantom{\Big|}$ \\
        \cline{2-4}
     &
    $-1$ &
    any &
    $\matRP^3\#\lens{t+3u}{u}\phantom{\Big|}$ \\
        \cline{2-4}
     &
    $-2$ &
    any &
    $\lens{5t+7u}{2t+3u}\phantom{\Big|}$ \\
        \cline{2-4}
     &
    $-1+\frac1n\phantom{\Big|}$ &
    $-1+\frac1m$ &
    $\lens{(2n+1)(2m+1)-4}{(2n+1)m-2}$\\
        \cline{3-4}
     &
    $\ne -2$ &
    $\ne -1,\ \ne -1+\frac1m\phantom{\Big|}\!\!$ &
    $\seiftre{S^2}21{2n+1}{-2}{t+u}u$\\
        \cline{2-4}
     &
    $-3$ &
    $\ne-1,\ \ne -1+\frac1m\phantom{\Big|}\!\!$ &
    $\seifdue{\matRP^2}21{t+u}{u}$ \\
        \cline{2-4}
     &
    $\begin{array}{c} \ne -1,\ \ne-3\\ \ne -1+\frac1n\end{array}$ &
    \!\!\!\!$\begin{array}{c} \ne -1,\ \ne-3\\
                        \ne -1+\frac1m \\
                \ne \frac{r}{s}\ {\rm if}\
                \frac{r}{s}\in\{1,-\frac{5}3\}\phantom{\Big|}
                    \end{array}$\!\!\!\!\!\! &
    $\begin{array}{l} \seifdue{D}21{r+s}s\bigcup\nolimits_{\tiny{\matr 0110}}
            \\ \seifdue{D}21{t+u}u\end{array}$ \\
        \hline
$-2$ &
    $-2$ &
    any &
    $\lens31\#\lens{t+2u}{u}\phantom{\Big|}$ \\ \cline{2-4}
&
    $-2+\frac1n$ &
    any &
    $\lens{3n(t+2u)-2t-u}{n(t+2u)-t-u}\phantom{\Big|}$ \\
        \cline{2-4}
&
    $\begin{array}{c} \ne-2\!\!\\ \ne-2+\frac1n\phantom{\Big|}\!\!\end{array}$ &
    $\begin{array}{c} \ne-2\!\!\\ \ne-2+\frac1m\phantom{\Big|}\!\!\end{array}$ &
    $\seiftre{S^2}32{r+2s}{-s}{t+2u}{-u}$ \\
        \hline
$-1$ &
    $-3$ &
    any &
    $\matRP^3\#\lens{t+3u}{u}\phantom{\Big|}$ \\
        \cline{2-4}
&
    $-3+\frac1n$ &
    any &
    $\lens{2n(t+3u)-t-u}{n(t+3u)+u}\phantom{\Big|}$ \\
        \cline{2-4}
&
    $\begin{array}{c} \ne-3\\ \ne-3+\frac1n\phantom{\Big|}\!\!\end{array}$ &
    $\begin{array}{c} \ne-3\\ \ne-3+\frac1m\phantom{\Big|}\!\!\end{array}$ &
    $\seiftre{S^2}21{r+3s}{-s}{t+3u}{-u}$ \\
        \hline
$0$ &
    $n$ &
    $-4-n$ &
    $\matRP^3\#\lens31\phantom{\Big|}$ \\
        \cline{3-4}
&
    &
    $-4-n+\frac1m$ &
    $\lens{6m-1}{2m-1}\phantom{\Big|}$ \\
        \cline{3-4}
&
    &
    $\begin{array}{c} \ne-4-n\\ \ne-4-n+\frac1m\phantom{\Big|}\!\!\end{array}$ &
    $\seiftre{S^2}2{-1}31{t+(n+4)u}{u}$ \\
        \cline{2-4}
&
    $\frac12+n$ &
    $-\frac{9}2-n$ &
    $\seifdueeul{\matRP^2}2131{-1}\phantom{\Big|}$ \\
        \cline{2-4}
&
    $\not\in\matZ$ &
    $\begin{array}{c}
          \not\in\matZ \\
          \neq -\frac rs -4 \ {\rm if}\ |s|=2
    \end{array}$ &
    $\begin{array}{l} \seifdue{D}{s}{r+2s}{u}{t+2u}\bigcup\nolimits_{\tiny{\matr
    01{-1}{-1}}}\phantom{\Big|}
        \\\seifdue{D}2131\phantom{\Big|}\end{array}$ \\
\end{tabular}
\nota{Classification of closed non-hyperbolic fillings of $N$}\label{closed:fill:table1}
\end{center}\end{table}

\begin{table}\begin{center}
\begin{tabular}{c|c||c||l}
$p/q$ & $r/s$ & $t/u$ & \qquad $M$ \\
        \hline\hline
$1$ & $1$ &
    $n$ &
    $T_{{\tiny{\matr{n+1}1{-1}0}}\phantom{\Big|}}\phantom{\Big|}$ \\
        \cline{3-4}
&
    &
    $\not\in\matZ$ &
    $\seifuno{A}{u}{t+u}\big/_{{\tiny{\matr 0110}}\phantom{\Big|}}\phantom{\Big|}$
    \\
        \hline
$-3/2$ & $-5/2$ &
    $-2$ &
    $\matRP^3\phantom{\Big|}$ \\
        \cline{3-4}
& &
    $-1$ &
    $\lens{13}{5}\phantom{\Big|}$ \\
        \cline{3-4}
& &
    $0$ &
    $\seifdueeul{\matRP^2}2131{-1}\phantom{\Big|}$ \\
        \cline{3-4}
& &
    $\begin{array}{c}-2+\frac1n\phantom{\Big|}\!\!\\ \ne-1\end{array}$ &
    $\seiftre{S^2}213{-1}{2n-1}2$ \\
        \cline{3-4}
& &
    $\begin{array}{c}\ne0\\ \ne -2\\ \ne-2+\frac1n\phantom{\Big|}\!\!\end{array}$ &
    $\begin{array}{l}\seifdue{D}2131\bigcup\nolimits_{\tiny{\matr 110{-1}}}
            \\ \seifdue{D}21{t+2u}u\end{array}$ \\
        \hline
$-4$ & $-1/2$ &
    $-1$ &
    $\lens{11}3\phantom{\Big|}$ \\
        \cline{3-4}
& &
    $-\frac12$ &
    $\seifdueeul{\matRP^2}2131{-1}\phantom{\Big|}$ \\
        \cline{3-4}
& &
    $0$ &
    $\lens{13}{5}\phantom{\Big|}$ \\
        \cline{3-4}
& &
    $\begin{array}{c}n\\ \ne-1\\ \ne0\end{array}$ &
    $\seiftre{S^2}213{-1}{2n+1}2$ \\
        \cline{3-4}
& &
    $\not\in\matZ\cup\{-\frac12\}$ &
    $\begin{array}{l} \seifdue{D}2131\bigcup\nolimits_{\tiny{\matr
    1110}}\phantom{\Big|}
            \\ \seifdue{D}21ut\phantom{\Big|}\end{array}$ \\
\end{tabular}
\nota{Classification of closed non-hyperbolic fillings of $N$ (continued)}\label{closed:fill:table2}
\end{center}\end{table}

\begin{table}\begin{center}
\begin{tabular}{c|c|c||l}
$p/q$ & $r/s$ & $t/u$ & \qquad $M$ \\
        \hline\hline
$-5$ & $-5$ & $-1/2$ &
    $\seifuno{A}21\big/_{{\tiny{\matr 0110}}\phantom{\Big|}}\phantom{\Big|}$ \\
    \hline
$-4$ & $-4$ & $-2/3$ &
    $\seifuno{A}21\big/_{{\tiny{\matr 1110}}\phantom{\Big|}}\phantom{\Big|}$ \\
    \hline
$-4$ & $-3/2$ & $-3/2$ &
    $T_{{\tiny{\matr {-3}1{-1}0}}\phantom{\Big|}}\phantom{\Big|}$ \\ \hline
$-4$ & $-1/3$ & $1$ &
    $\seifdue{D}2121
        \bigcup\nolimits_{{\tiny{\matr 01{-1}{-1}}}\phantom{\Big|}\!\!}
        \seifdue{D}2131\phantom{\Big|}$ \\ \hline
$-8/3$ & $-3/2$ & $-3/2$ &
    $\seifdue{D}2121
        \bigcup\nolimits_{{\tiny{\matr 01{-1}{-1}}}\phantom{\Big|}\!\!}
        \seifdue{D}2131\phantom{\Big|}$ \\ \hline
$-5/2$ & $-5/2$ & $-4/3$ &
    $\seifuno{A}21\big/_{{\tiny{\matr 2110}}\phantom{\Big|}}\phantom{\Big|}$ \\
    \hline
$-5/2$ & $-5/3$ & $-5/3$ &
    $\seifdue{D}2121
        \bigcup\nolimits_{{\tiny{\matr {-1}10{-1}}}\phantom{\Big|}\!\!}
        \seifdue{D}2131\phantom{\Big|}$ \\ \hline
$-7/3$ & $-7/3$ & $-3/2$ &
    $\seifuno{A}21\big/_{{\tiny{\matr 1110}}\phantom{\Big|}}\phantom{\Big|}$ \\
    \hline
$1$ & $2$ & $2$ &
    $\seiftreeul{S^2}213171{-1}\phantom{\Big|}$ \\ \hline
$1$ & $2$ & $3$ &
    $\seiftreeul{S^2}214151{-1}\phantom{\Big|}$ \\ \hline
$1$ & $2$ & $4$ &
    $\seiftreeul{S^2}313141{-1}\phantom{\Big|}$ \\ \hline
$1$ & $2$ & $5$ &
    $\seifdue{D}2121
        \bigcup\nolimits_{{\tiny{\matr 0110}}\phantom{\Big|}\!\!}
        \seifdue{D}2131\phantom{\Big|}$ \\ \hline
$1$ & $3$ & $3$ &
    $\seifdue{D}2121
        \bigcup\nolimits_{{\tiny{\matr 120{-1}}}\phantom{\Big|}\!\!}
        \seifdue{D}2131\phantom{\Big|}$ \\ \hline
$2$ & $2$ & $2$ &
    $\seifdue{D}2121
        \bigcup\nolimits_{{\tiny{\matr 23{-1}{-2}}}\phantom{\Big|}\!\!}
        \seifdue{D}2131\phantom{\Big|}$ \\
\end{tabular}
\nota{Classification of closed non-hyperbolic fillings of $N$ (continued)}\label{closed:fill:table3}
\end{center}\end{table}

The next remark is important to understand the statement of
Theorem~\ref{three:cusps:teo}.

\begin{rem}\label{three:cusps:explain:rem}
\emph{We have compiled Tables~\ref{closed:fill:table1} to~\ref{closed:fill:table3} trying to keep them concise but
easily usable. In particular:
\begin{itemize}
\item we have always excluded the case where one of the filling
coefficients is $\infty$, because the first point of Theorem~\ref{three:cusps:teo}
is already sufficient in this case to determine which lens space the result is;
\item whenever a certain value $p/q\in\{-3,-2,-1,0\}$ is chosen, the corresponding
set of lines in Table~\ref{closed:fill:table1} contains a complete classification of the resulting $M$, but
one may have to permute $r/s$ and $t/u$ to get the result;
\item the tables
always provide the decomposition of $M$ into prime connected
summands and the JSJ decomposition of these summands, except when a connected
sum involves a lens space and the parameters are such that the lens
space is $S^3$;
\item the triples listed in Table~\ref{closed:fill:table3} do not appear
in Tables~\ref{closed:fill:table1} and~\ref{closed:fill:table2};
\item there is some overlapping
between the cases covered in Tables~\ref{closed:fill:table1} and~\ref{closed:fill:table2}. The
recognition of $M$ is of course coherent (note for instance that $N(1,1,n) = T_{\tiny\matr{n+1}1{-1}0}$
is also a Seifert manifold for $n\in\{-3,-2,-1,0\}$, see Lemma~\ref{graph:lem});
\item in Table 3 the
elements of each triple are arranged in increasing order, and the triples themselves are
arranged in a lexicographically increasing order.
\end{itemize}}
\end{rem}

\paragraph{Sketch of the proofs}
The proof of Theorems~\ref{one:cusp:teo} to \ref{three:cusps:teo} goes as follows: we first prove
Theorem~\ref{one:cusp:teo} using spines. Actually, we prove a more refined version of it, in which we
determine the image in $N(\alpha)$ of the bases of the non-filled boundary tori. Using spines, we
also recognize $N(1,1),$ $N(2,2,2),$ and $N(1,3,3)$. We then use Proposition~\ref{repetitions:prop} below
(discussed in Section~\ref{proofs:section})
to recognize all the other fillings
listed in Tables~\ref{two:cusps:table} to~\ref{closed:fill:table3}.
Next, we use the Gromov-Thurston $2\pi$ theorem to show
that if $M$ is a filling of $N$, but not one of those already recognized to
be non-hyperbolic, then $M$ admits a complete metric with
negative sectional curvature.
An application of the orbifold theorem~\cite{Bo-Po} then implies that
$M$ is actually hyperbolic. This is because our chain-link
complement $N$ admits an orientation-preserving order-$2$ involution, given by a rotation
of angle $\pi$ around a horizontal circle in Fig.~\ref{chain:fig}.
Its effect on the homology of each cusp is given by multiplication by $-1$,
therefore the involution extends to every Dehn filling of $N$.
The orbifold theorem then guarantees that each filling
of $N$ satisfies the geometrization conjecture.

\begin{prop}\label{repetitions:prop}
$N(-3/2) = N(-4)$ and
the following holds for all $\alpha,\beta\in\matQ\cup\{\infty\}$:
\begin{eqnarray}
N\left(-\frac 32,\alpha,\beta\right) & = & N\left(-4,-\frac{\alpha+1}{\alpha+2},-\beta-3\right)\label{511right:eqn} \\
N\left(-\frac 52,\alpha, \beta\right) & = & N\left(-\frac 52,-\frac {\alpha+3}{\alpha+2},-\frac {2\beta+3}{\beta+1}\right)\label{512:eqn} \\
N\left(-\frac 32,\alpha, \beta\right) & = & N\left(-\frac 32,-\frac {2\alpha+5}{\alpha+2},-\frac{2\beta+5}{\beta+2}\right)\label{511left:eqn} \\
N\left(-\frac 12,\alpha, \beta\right) & = & N\left(-\frac 12,-\alpha -4,-\beta-4\right)\label{513:eqn} \\
N\left(1,2,\alpha\right) & = & N\left(1,2,-\alpha+2\right)\label{514:eqn}   \\
N\left(1,-4,\alpha\right) & = & N\left(1,-4,\frac 1\alpha\right) \label{515:eqn}
\end{eqnarray}
Relations~(\ref{512:eqn}) to~(\ref{515:eqn}) are induced by automorphisms of
$N(-5/2)$, $N(-3/2)$, $N(-1/2)$, $N(1,2)$, and $N(1,-4)$, respectively.
Moreover~(\ref{512:eqn}) has order $3$, while~(\ref{511left:eqn})-(\ref{515:eqn}) have order $2$.
\end{prop}
We note that $N$ and
the manifolds involved in
Proposition~\ref{repetitions:prop} are well-known hyperbolic ones.
Using SnapPea one sees that $N$ is the manifold called the
``magic'' one in~\cite{GoWu-TA} and listed as $M6_1^3$
in~\cite{Ca-Hi-We}.

\begin{wrapfigure}[12]{R}{2in}
  \centering
  \vspace{.5 cm}
  \mettifig{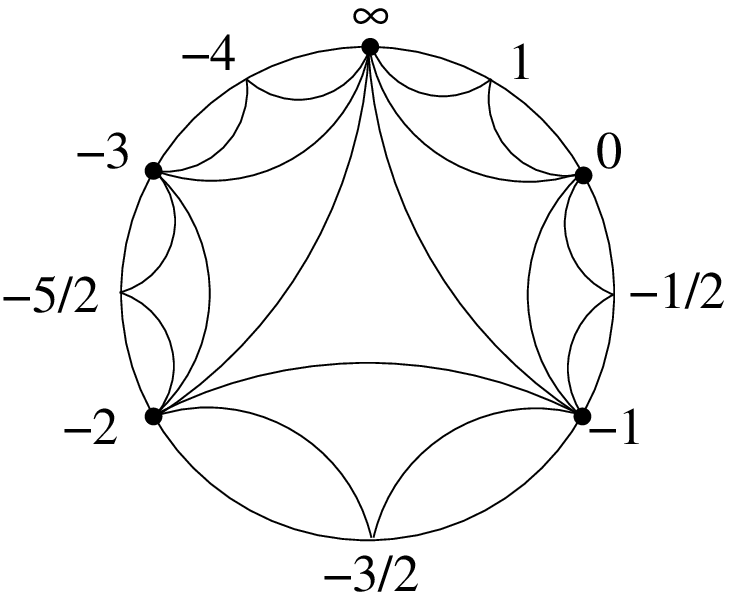, width = 4.5 cm}
  \nota{Important fillings of $N$.}
  \label{N:fig}
\end{wrapfigure}
Some important fillings of $N$ are link
complements in the 3-sphere:
$N(-4) = N(-3/2)$ of the Whitehead
sister link,
$N(-5/2)$ of the Berge link, $N(-1/2)$ of the
$3/10$-bridge link, $N(1)$ of the Whitehead link. In addition
$N(1,2)$ is the figure-8 knot complement and $N(1,-4)$ is the
figure-8 knot sister. See the Appendix for more details on these
manifolds and their fillings. We show the
slopes involved in these fillings in Fig.~\ref{N:fig} where we
arrange them on $\partial\matH^2=\matQ\cup\{\infty\}$ using the
disc model of $\matH^2$, and we connect by lines slopes having
distance $1$, thus showing a partial picture of the so-called
Farey tessellation.

\subsection{Some consequences} \label{consequences:subsection}

We illustrate here some facts that can be deduced from our main
statements.

\paragraph{Seifert fillings}
We begin by proving Theorem~\ref{Seifert:teo} as a
consequence of the following:
\begin{prop} \label{Seifert:prop}
For every integer $k\neq 0$, among the manifolds of the form
$N\left(1+\frac 1h, \frac{kh+1}{k(h+1)+1}\right)$, with
$h\in\matZ$, there are infinitely many distinct
hyperbolic ones whose exceptional slopes are precisely
$\infty,-3,-2,-1,0$. Moreover
\begin{itemize}
\item The $\infty$-filling is $S^3$; \item The $(-3)$-filling and
the $0$-filling are toroidal; \item The $(-2)$-filling and the
$(-1)$-filling are small Seifert manifolds with first homology
respectively $\matZ/_{(k+3)\matZ}$ and $\matZ/_{(k+2)\matZ}$.
\end{itemize}
\end{prop}

\begin{proof}
A manifold $M$ as described is hyperbolic by
Theorem~\ref{two:cusps:teo}, except for finitely many values of $h$.
Moreover as $h$ goes to $\infty$ the
volume of $M$ tends to the volume of $N$ from below, so there are
infinitely many different $M$'s. For the other assertions, setting
$r/s = (h+1)/h$ and $t/u = (kh+1)/(k(h+1)+1)$, we have
$$tr-us=\det \matr rust = \det \matr {h+1}{1+k(h+1)}h{1+kh} = \det \matr
{h+1}1h1 = 1$$
which implies (by the first assertion of
Theorem~\ref{three:cusps:teo}) that $M(\infty) = S^3$. As $h$
tends to $\infty$ (while $k$ stays fixed) both $r/s$ and $t/u$
tend to $1$, and it easily follows from
Table~\ref{closed:fill:table1} that only finitely many values of
$h$ can yield manifolds of non-generic type for $M(-3), M(-2),
M(-1)$, and $M(0)$. Then for infinitely many $h$'s the topology of
the fillings is as required. We are left to prove the assertion
concerning homology.

It is easy to see that the number $\nu = |q_1p_2p_3 + q_2p_3p_1 +
q_3p_1p_2|$ equals the cardinality of $H_1\big(S^2, (p_1,q_1),
(p_2,q_2), (p_3,q_3)\big)$ except when $\nu=0$, in which case the
group is infinite. When the Seifert manifold is obtained from
surgery on a knot in $S^3$, its first homology group must be
$\matZ/_{\nu\matZ}$. Therefore:
\begin{eqnarray*}
H_1\big(M(-2)\big) & = & H_1\big(S^2, (3,2), (r+2s, -s), (t+2u,
-u)\big) \\
& = & \matZ/_{(2(r+2s)(t+2u) - 3s(t+2u) -3u(r+2s))\matZ} \\
& = & \matZ/_{(2rt + ru + st - 4su)\matZ}  = \matZ/_{((r-s)(u-t) + 3(rt-su))\matZ} = \matZ/_{(k+3)\matZ}
\end{eqnarray*}
and similarly $H_1\big(M(-1)\big) = \matZ_{(k+2)\matZ}$.
\end{proof}

\vspace{2pt}\noindent\textit{Proof of Theorem~\ref{Seifert:teo}}.\
Fix $n$, assume first $n\neq -4$, and apply
Proposition~\ref{Seifert:prop} with $k=-n-4$. Each of the
manifolds given by the
proposition is the complement of a knot having 4 consecutive
integer exceptional slopes, say $a,a+1,a+2,a+3$. Since the first
homology of a $b$-surgery on a knot is $\matZ/_{b\matZ}$, we
deduce that
$$\matZ/_{(a+1)\matZ}=\matZ/_{(k+3)\matZ},\qquad
\matZ/_{(a+2)\matZ}=\matZ/_{(k+2)\matZ}.$$ Now $k+3=-(n+1)$ and
$k+2=-(n+2)$, whence $a+1=\pm(n+1)$ and $a+2=\pm(n+2)$, whose only
solution is $a=n$. This implies that we have infinitely many knots
with exceptional slopes $n,n+1,n+2,n+3$
as required. For $n=-4$ it is sufficient to take the knots for
$n=1$ and consider their mirror images.\finedimo

\paragraph{Knots with homeomorphic complements}
We discuss here Theorems~\ref{knots1:teo} and~\ref{knots2:teo}. As
usual, we translate them into a statement about Dehn fillings: two
inequivalent knots in $M$ having the same $M'$ as their
complement
correspond to inequivalent slopes in $\partial M'$ both giving $M$
via Dehn filling. Two slopes are \emph{equivalent} if there is an
automorphism of $M'$ sending one to the other.

We also recall some other definitions from~\cite{cosmetic}. Two
slopes $\alpha$ and $\beta$ such that $M'(\alpha)$ is homeomorphic
to $M'(\beta)$ are called \emph{cosmetic}. Taking on $M'(\alpha)$ and $M'(\beta)$ an orientation
that extends one and the same orientation of $M'$, we say that
$\alpha$ and $\beta$ are \emph{truly cosmetic} if there is an
orientation-preserving homeomorphism from $M'(\alpha)$ to
$M'(\beta)$, and \emph{reflectively cosmetic} if there is an
orientation-reversing one (both possibilities can hold).

The only known $1$-cusped hyperbolic $M$ with a pair of
inequivalent cosmetic slopes was found in~\cite{cosmetic} and it
is $M=N(-5/2,-2/5)$ which admits the cosmetic pair of slopes
$-2,-1$ giving $\lens{49}{18}$, as one sees from
Table~\ref{closed:fill:table1} or Table~\ref{two:cusps:table3} in
the Appendix. The two slopes are not truly cosmetic. In fact, it
is conjectured~\cite{cosmetic, Kirby} that inequivalent slopes on
a $1$-cusped $M$ are never truly cosmetic.

Cosmetic phenomena seem to be more frequent on manifolds with more
cusps. All the $2$-cusped $M$'s with distinct slopes
$\alpha,\beta$ on a torus $T\subset\partial M$ such that
$M(\alpha)=M(\beta)$ is a solid torus are described by Berge and
Gabai in~\cite{Be, Ga2}. Infinitely many such pairs of slopes are
inequivalent. This also follows from the next result, proved in
Subsection~\ref{number:subsection} of the Appendix, which
implies Theorem~\ref{knots1:teo} and~\ref{knots2:teo} too.

\begin{teo}~\label{cosmetic:teo}
Let $M$ be $N$ or a hyperbolic $N(\gamma)$, and let $\alpha,\beta$
be a cosmetic pair of inequivalent slopes on $M$. Then
$M,\alpha,\beta$ are as in Table~\ref{cosmetic:table}.
\end{teo}
\begin{table}
\begin{center}
\begin{tabular}{c|c|c|c|c|c}
$M$ & $\alpha$ & $\beta$ & $M(\alpha)$ & truly & reflectively \\
\hline\hline
$N$ & $-4$ & $-3/2$ & Whitehead sister & no & yes \phantom{\Big|}\\
\hline
$\begin{array}{c} N(-3+1/n) \\ n\in\matZ\setminus\{0,\pm 1, 2\} \end{array}$ &
$\infty$ & $-1$ & $D\times S^1$ & yes & yes \phantom{\Big|} \\
\hline
$\begin{array}{c} N(-2+1/n) \\ n\in\matZ\setminus\{0,\pm 1, \pm 2\} \end{array}$ &
$\infty$ & $-2$ & $D\times S^1$ & yes & yes \phantom{\Big|} \\
\hline
$N(-12/5)$ & $-2$ & $-1$ & $\seifdue D2131$ & yes & no \phantom{\Big|} \\
\hline
$N(-6)$ & $-1$ & $0$ & $\seifdue D2131$ & no & yes \phantom{\Big|} \\
\hline
$N(-4/3)$ & $-3$ & $-1$ & $\seifdue D2152$ & no & yes \phantom{\Big|}
\end{tabular}
\nota{Cosmetic slopes on $N$ and $N(\gamma)$}\label{cosmetic:table}
\end{center}
\end{table}

Using SnapPea, one sees that the two cosmetic knots in the
Whitehead link sister corresponding to the first line of
Table~\ref{cosmetic:table} are the two shortest geodesics, with
lengths $0.96242\ldots$ and $1.06128\ldots$.

\paragraph{Non-hyperbolic closed fillings}
The following result is easily proved
as a consequence of Theorem~\ref{three:cusps:teo} and Lemma~\ref{graph:lem}:

\begin{cor} \label{all:manifolds:cor}
The set of all non-hyperbolic closed manifolds obtained
by Dehn surgery on the chain-link consists of:
\begin{itemize}
\item $S^3$, $S^2\times S^1$, $\lens pq$,
\item $\matRP^3\#\lens pq$, $\lens31\#\lens pq$,
$\matRP^3\#(S^2\times S^1)$, $\lens31\#(S^2\times S^1)$,
\item $\seiftreeul {S^2}21pqrsb, \seiftreeul {S^2}31pqrsb, \seifdueeul {\matRP^2}21pqb$,
\item $\seifdue D21pq \bigcup_X \seifdue D21rs$,
\item $\seifdue D2131 \bigcup\nolimits_{{\tiny{ \matr 11n{n\pm 1} }}} \seifdue Dpqrs$,
\item $T_{\tiny{\matr n1{-1}0}}$,
\item $\seifuno Apq /_Y$,
\item $\seifdue{D}2121\bigcup\nolimits_{{\tiny{\matr 120{-1}}}}
        \seifdue{D}2131$,
\item $\seifdue{D}2121\bigcup\nolimits_{{\tiny{\matr 23{-1}{-2}}}}
        \seifdue{D}2131$,
\end{itemize}
where $(p,q)$ and $(r,s)$ are any coprime pairs with $p>q>0$ and $r>s>0$, $b,n$ are any
integers, and
$X,Y\in\GL_2(\matZ)$
are any matrices whose top-right entry equals $\pm1$, with ${\rm det}(Y)=-1$.
\end{cor}

Concerning this corollary,  note that the top-right entry of a matrix $X$ or $Y$ as in
the statement is the geometric intersection number of the fibres of the Seifert fibrations that are glued
together. The only graph manifolds of the list for which this intersection number is not $\pm 1$ are
$N(2,2,2)$ and $N(1,3,3)$. We also note that  not all the manifolds listed can be obtained as fillings of a
\emph{hyperbolic}  1-cusped $N(\alpha, \beta)$. For instance, the only reducible manifolds thus obtained
are $S^2\times S^1$ (with  $\alpha = 1/\beta$ and $\alpha,\beta\in\matQ\setminus\{-3,-2,-1,1\}$) and
$\matRP^3\#\lens 31$  (with $\alpha = n$, $\beta = -4-n$ and $n\in\matZ\setminus\{-4,-3,-2,-1,0\}$).

\paragraph{Chirality}
The following result confirms the experimental fact~\cite{Ca-Hi-We} that amphichiral manifolds are quite
sporadic, at least among those with small volume.
\begin{prop} \label{amphichiral:prop}
The only amphichiral hyperbolic manifolds with boundary obtained
as fillings of $N$ are $N(-1/2)$, $N(1,2)$, and
$N(1,-4)$.
\end{prop}
The three manifolds are the figure-$8$ knot complement, its sister,
and the $3/10$-bridge link complement, see Fig.~\ref{N:fig} above.
The proof is contained in Subsection~\ref{number:subsection} of the Appendix.

\section{Proofs of the main results}\label{proofs:section}
We give here a unified proof of Theorems~\ref{one:cusp:teo},
\ref{two:cusps:teo}, and \ref{three:cusps:teo}. As already
mentioned, the ingredients we use are
Proposition~\ref{repetitions:prop}, a certain refined version of
Theorem~\ref{one:cusp:teo} (proved below as
Proposition~\ref{recognized:fillings:prop}), and then an elaborate
argument which uses the Gromov-Thurston $2\pi$-theorem and SnapPea
for a finite enumeration.

\paragraph{Hyperbolic fillings with symmetries}
One can see using SnapPea (and check by hand) that
the partial Dehn fillings
$$N(-3/2)\qquad N(-4)\qquad N(-5/2)\qquad N(-1/2)\qquad N(1,2)\qquad N(1,-4)$$
are all hyperbolic, that the first two of them are isometric to each other,
and that in all cases the symmetry group contains a non-trivial cyclic
subgroup that leaves invariant each boundary component. Using the precise form
of these non-trivial isometries and their action
on the homology bases on the boundary, one can
with some patience deduce the identities of
Proposition~\ref{repetitions:prop}.

\paragraph{Identities between graph manifolds}
We slightly relax in this section our notation
$\big(F,({\scriptstyle p_1},{\scriptstyle q_1}),\ldots,
({\scriptstyle p_k},{\scriptstyle q_k}),b\big)$
for Seifert manifolds, allowing also $|p_i|\leqslant 1$,
which can lead to lens spaces or to connected sums of lens
spaces and/or Seifert manifolds. We will always use the conventions
of Section~\ref{main:section} on homology bases to encode gluings.
The following lemma, which collects well-known~\cite{FoMa} or easy
properties of graph manifolds, is repeatedly (and often tacitly) used below.

\begin{lemma} \label{graph:lem}
Suppose $F$ is a surface, $(i,j)$ and $(i_n,j_n)$ are coprime pairs,
$k\in\matZ$, $X$ is any graph manifold, and $\tiny \matr abcd\in\GL_2(\matZ)$.
Let $(i',j')$ be such that $ij'-ji'=\pm 1$ and, when $\partial F\neq\emptyset$,
denote by $F'$ the surface obtained from $F$ by capping off one boundary circle.
Then the following equalities hold:
\begin{description}
\item[(Seifert manifolds)]
\begin{eqnarray}
\big(F, ({\scriptstyle i_1},{\scriptstyle j_1}),({\scriptstyle i_2},{\scriptstyle j_2}),\ldots\big)
    & = &\big(F, ({\scriptstyle i_1},{\scriptstyle -j_1}),({\scriptstyle i_2},{\scriptstyle -j_2}),\ldots\big)
    \label{minus:eqn} \\
\seifdueeul F{i_1}{j_1}{i_1}{j_2}\ldots & = & \seifdueeul F{i_1}{j_1+k i_1}{i_2}{j_2-k i_2}\ldots \label{from:D2:to:D2:eqn}\\
\big(F, ({\scriptstyle i_1},{\scriptstyle j_1}),({\scriptstyle i_2},{\scriptstyle j_2}),\ldots\big)
    & = &\big(F, ({\scriptstyle i_1},{\scriptstyle j_1+k i_1}),({\scriptstyle i_2},{\scriptstyle j_2}),\ldots\big)
    \quad {\rm if\ }\partial F \neq \emptyset \label{from:F:to:F:eqn}\\
\seifdue {S^2}{i}{j}{i_1}{j_1} & = & \lens {i_1j+j_1i}{i_1j'+j_1i'}
    \label{from:D1:to:lens:eqn}\\
\seifdueeul F1{k}{i}{j}\ldots & = & \seifunoeul {F}{i}{j+k i}\ldots \label{from:D2:to:D1:eqn}
\end{eqnarray}

\item[(Graph manifolds)]
\begin{eqnarray}
\seifunoeul F{i}{j}\ldots\bigu abcd X & = & \seifunoeul F{i}{j+k i}\ldots\bigu
    {a+k b}b{c+k d}d X \label{from:D2X:to:D2X:eqn}\\
X\bigu abcd \seifunoeul F{i}{j}\ldots & = & X\bigu ab{c-k a}{d-k b}
    \seifunoeul Fi{j+ki}{\ldots} \label{from:XD2:to:XD2:eqn}\\
\seifuno D{i}{j}\bigu abcd \seifeul F\ldots & = &
    \seifunoeul {F'}{ai-bj}{ci-dj}\ldots \label{from:D1D2:to:D2:eqn}
\end{eqnarray}

\item[(Double fibrations)]
\begin{eqnarray}
\seifdue D2121 \bigu abcd  X & = &
(S\timtil S^1)\bigu b{b-a}d{d-c} X
    \label{from:D2X:to:SX:eqn} \\
X\bigu abcd \seifdue D2121 & = & X \bigu{a+c}{b+d}{-a}{-b} S\timtil S^1
    \label{from:XD2:to:XS:eqn}
\end{eqnarray}

\item[(Reducible manifolds)]
\begin{eqnarray}
\seiftre {S^2}{i_1}{j_1}{i_2}{j_2}01 & = & \lens {i_1}{j_1}\# \lens {i_2}{j_2}
    \label{from:D2:to:lens:plus:lens:eqn}  \\
\seifdue D01{i}{j}\bigu abcd \seifeul F\ldots & = & \lens {i}{j}\#\seifunoeul {F'}bd\ldots
    \label{from:DX:to:lensX:eqn}
\end{eqnarray}

\item[(Torus bundles)]
\begin{eqnarray}
\big(A\times S^1\big) \bigb abcd & = & T \bigb {-a}{-b}cd
\label{from:A:to:T:eqn} \\
\seiftreeul {S^2}213161{-1} & = & T \bigb 11{-1}0
\label{from:Seif:to:T:eqn1} \\
\seiftreeul {S^2}214141{-1} & = & T \bigb 01{-1}0
\label{from:Seif:to:T:eqn2} \\
\seiftreeul {S^2}313131{-1} & = & T \bigb {-1}1{-1}0
\label{from:Seif:to:T:eqn3} \\
\seifeul K1 & = & T \bigb {-2}1{-1}0. \label{from:Seif:to:T:eqn4}
\end{eqnarray}
\end{description}
The homeomorphisms proving these equalities respect the homology bases on the
boundary, except in formulae~(\ref{minus:eqn}) and~(\ref{from:F:to:F:eqn}).
Relations (\ref{from:D2X:to:D2X:eqn}) and (\ref{from:XD2:to:XD2:eqn})
are also valid if $\tiny \matr abcd$ is used to glue together
two boundary components of the same Seifert block $\seifunoeul F{i}{j}\ldots$.
\end{lemma}

\paragraph{Basic recognition result}
Recall that we are always denoting by $N$ the complement
of the chain-link with 3 components. We establish now
the following:

\begin{prop} \label{recognized:fillings:prop}
The following homeomorphisms hold:
\begin{eqnarray}
  N\left(\infty,\frac rs,\frac tu\right) & = & \seifdue {S^2}rs{-u}t\label{521:eqn} \\
  N\left(-3,\frac rs,\frac tu\right) & = & \seifdue D21{r+s}s\bigu 0110 \seifdue D21{t+u}u\label{522:eqn} \\
  N\left(-2,\frac rs,\frac tu\right) & = & \seiftre {S^2}32{r+2s}{-s}{t+2u}{-u} \label{523:eqn} \\
  N\left(-1,\frac rs,\frac tu\right) & = & \seiftre {S^2}21{r+3s}{-s}{t+3u}{-u} \label{524:eqn}
\end{eqnarray}
\begin{eqnarray}
  N\left(0,\frac rs,\frac tu\right)  & = & \seifdue
  Ds{r+2s}u{t+2u}\bigu 01{-1}{-1} \seifdue D2131 \label{525:eqn} \\
  N\left(1,1,\frac tu\right)         & = & \seifuno Au{t+u}\bigb 0110
  \label{526:eqn} \\
  N\left(1,3,3\right)                & = & \seifdue{D}2121 \bigu 120{-1}    \seifdue{D}2131 \label{527:eqn} \\
  N\left(2,2,2\right)                & = & \seifdue{D}2121 \bigu 23{-1}{-2} \seifdue{D}2131 \label{528:eqn}
\end{eqnarray}
Equalities (\ref{521:eqn}) to (\ref{526:eqn}) also hold for partial fillings:
if a symbol $r/s$ or $t/u$ is missing on the left-hand side (so
the relative cusp is not filled),
an open neighborhood of the corresponding fibre on the
right-hand side should be removed.
\end{prop}

\begin{proof}
  The first equality is easy: doing $\infty$-surgery on a component
  of a link means canceling this component. So $N(\infty)$ is the
  complement of the Hopf link, which is homeomorphic to $T\times I$,
  and it is easy to see that the actual homeomorphism leads to
  relation~(\ref{521:eqn}).

  The other equalities are proved using spines.
  Recall that a \emph{spine} of a PL manifold $M$ is a subpolyhedron onto which
  $M$ collapses. Moreover, if $M$ is cusped, the 2-skeleton of the cellularization
  dual to an ideal triangulation of $M$ is always a spine $P$ of $M$.
  The \emph{singular} (non-surface) points of $P$ give a 4-valent
  graph $S(P)$ such that $P\setminus S(P)$ consists of discs (called the \emph{faces} of $P$).
  A spine of this sort is called \emph{standard}, and it is
  known to determine $M$. Moreover, $P$ itself can be reconstructed
  by a neighbourhood of $S(P)$ in $P$.
  Below we will employ also many non-standard spines.

The usual triangulation technique for complements of alternating links
shows that our chain-link complement $N$ is obtained by pairing the faces
of two triangular prisms. Subdividing the prism and dualizing we see
that $N$ has a standard spine with $6$ vertices. One can actually
show that $6$ is the minimal number of vertices for a spine of $N$,
and that there are several
spines with $6$ vertices, related to the original one by the
moves of the Matveev-Piergallini calculus. One such spine, denoted by $P$ and
shown in Fig.~\ref{P:fig}, is particularly suitable for the recognition of fillings.
We could give a constructive proof of the fact that $P$ is a spine of $N$,
but we confine ourselves to a quick indirect proof.

  \begin{figure}
    \begin{center}
      \mettifig{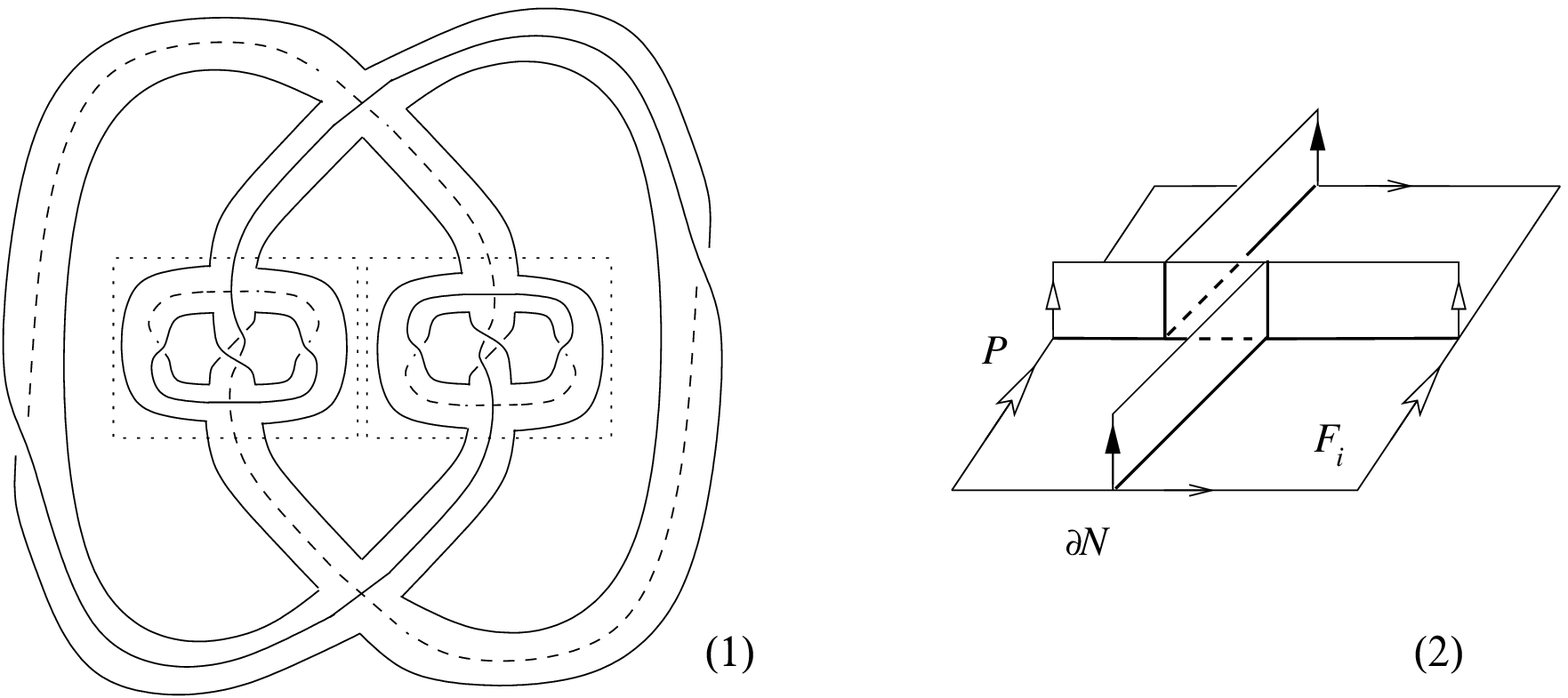, width = 13cm}
      \nota{A neighbourhood of $S(P)$ and
      a neighbourhood of $\overline{F_i}$ (edges with the same label are identified)}
      \label{P:fig}
    \end{center}
  \end{figure}

\begin{lemma}
Let $P$ be the standard polyhedron with a neighbourhood of $S(P)$ as in
Fig.~\ref{P:fig}-(1). Then:
\begin{enumerate}
\item $P$ is a spine of $N$;
\item $P$ has two faces $F_1$ and $F_2$ such that $\overline{F_i}$ is a torus parallel
to a component $T_i$ of $\partial N$, and the rest of $P$ is attached to
$\overline{F_i}$ along a $\theta$-shaped graph as shown in Fig.~\ref{P:fig}-(2);
the attachment takes place on the side of $\overline{F_i}$ opposite
to the parallelism with $T_i$;
\item There are $4$ other faces in $P$; a hexagon $E$ and a decagon $C$,
each incident to both $\overline{F_i}$'s along an edge,
and two squares $Q_1$ and $Q_2$, with $Q_i$ incident to $\overline{F_i}$ along an edge.
\end{enumerate}
\end{lemma}

\begin{proof}
Points (2) and (3) are proved by direct inspection.
Let $M$ be the manifold $P$ is a spine of. By point (2) $M$ has
at least 3 boundary components, and dual to $P$ there is
an ideal triangulation of $M$ involving 6 tetrahedra.
We now claim $M$ is (complete, finite-volume) hyperbolic. This
is sufficient to prove that $M=N$, because $N$ is the unique
hyperbolic manifold with at least three cusps which can be
triangulated with at most 6 tetrahedra~\cite{Ca-Hi-We}.

We prove the claim recalling that Thurston has associated certain
hyperbolicity equations to each triangulated manifold, and that
Casson has shown
that a solution of the angle part
of the consistency equations is already sufficient to ensure
hyperbolicity~\cite{lackenby}.
In terms of spines, we must assign angles to the corners of the faces,
and we do this as follows:
(1) We put $\pi/3$ on the corners of $E$, and on the 6 corners opposite to those of $E$;
(2) We put $\pi/2$ on the corners of the $Q_i$'s and on those opposite to them;
(3) We put $\pi/6$ on all other corners. One readily sees that
angles at opposite corners are always the same and the sum is $2\pi$ on all faces,
so we have solved the angle equations.
\end{proof}

Let us resume the proof of Proposition~\ref{recognized:fillings:prop}.

\paragraph{Homology bases}
We now isotope $\overline{F_i}$ so that it coincides with $T_i\subset \partial N$,
and note that the complement in $T_i$ of the $\theta$-graph
$\partial F_i$ is a disc, so $\partial F_i$ contains three
unoriented loops giving three slopes on $T_i$
with mutual distance $1$. We can therefore choose a basis
of $H_1(T_i)$ so that $\partial F_i$ contains the slopes
$-1$, $-2$, and $\infty$. More precisely, we arrange so that
$-1$ (respectively, $-2$ and $\infty$) is the slope
obtained by erasing from $\partial F_i$ the edge which belongs
to $Q_i$ (respectively, $E$ and $C$). We will also say that
$-1$ (respectively, $-2$ and $\infty$) is the slope
associated to $Q_i$ (respectively, $E$ and $C$).

We will now proceed using these homology bases,
recognize many fillings of $N$, and only later show that
these bases are actually (up to sign) the natural ones.
For this reason we recognize again $N(\infty)$, because
it is not obvious at this stage that $\infty$ represents
the same slope it did with respect to the meridian-longitude basis.

\paragraph{General recognition strategy}
Our recognition of $N(\alpha)$ for $\alpha\in\{\infty,-3,-2,-1,$ $0\}$
is based on the construction of a spine $P(\alpha)$.
We first illustrate here how to build $P(\alpha)$ for arbitrary $\alpha$.
The same method also gives a spine of $N(\alpha,\beta)$,
that we will use to recognize $N(1,1)$,
but it does not apply to the closed fillings $N(\alpha,\beta,\gamma)$.
To recognize $N(1,3,3)$ and $N(2,2,2)$ we have
constructed spines of the punctured manifolds via a slightly different
technique, sketched below.

  The spine $P(\alpha)$ of $N(\alpha)$ is constructed from $P$ by
  removing the (open) hexagon $F_1$, thus getting
  a $\theta$-shaped ``boundary,'' and attaching a polyhedron
  $Q(\alpha)$ having $\theta$-shaped boundary too, with
  a map matching the two $\theta$'s. More precisely:
  consider $N(\alpha) = N\cup V$ where $V$ is an open solid torus,
    glued to $N$ along $T_1\subset\partial N$.
Since $\overline{F_1}=T_1$, then
  $N(\alpha)\setminus P$ is the union of $V$ and two copies of $T\times (0,1]$'s.
  Suppose we have a polyhedron $Q(\alpha)\subset \overline V$ so that
  (i) $Q(\alpha)\cap \partial V$
  is contained in the $\theta$-shaped graph $\partial F_1$,
    (ii) $V\setminus Q(\alpha)$ is an open ball.
  By (i), the complement of $P\cup Q(\alpha)$ inside $N(\alpha)$
  is the union of two $T\times (0,1]$'s and one open ball.
    By (ii), $F_1$ is again an (open) face of $P\cup Q(\alpha)$, which separates the ball
  from one $T\times (0,1]$. Therefore $(P\cup Q(\alpha))\setminus F_1$
  is a spine of $N(\alpha)$, which we call $P(\alpha)$.

  We explain now how to construct $Q(\alpha)$ for all $\alpha$
  (even if we will only be using the construction for a few $\alpha$'s).
  If $\alpha\in\{-2,-1,\infty\}$
  then the meridian of $V$ is isotopic to one of the three loops in
  $\partial F_1$, and
  we are in a ``degenerate'' case, because $Q(\alpha)$ can be taken to be just a meridinal
  disc of $V$, as in Fig.~\ref{mobius_triplet:fig}-(1), so $P(\alpha)$ is not standard.
  If $\alpha\in\{-3,-3/2,0\}$, the $\theta$-shaped graph $\partial F_1$ does not contain
  the meridian of $V$, but it intersects the meridian in two points.
  We can then take as $Q(\alpha)$
  the M\"obius triplet shown in Fig.~\ref{mobius_triplet:fig}-(2).
  For all remaining $\alpha$'s the $\theta$-shaped graph intersects the meridian in more than two
  points. Therefore we change the $\theta$-shaped
  boundary of the M\"obius triplet via some \emph{flips}
  (see Fig.~\ref{mobius_triplet:fig}-(3,4)) in order to get a $\theta$-shaped graph
  which matches $\partial F_1$, each flip
  adding a vertex to $Q(\alpha)$ as shown in Fig.~\ref{mobius_triplet:fig}-(5).
  (This $Q(\alpha)$ is dual to a layered triangulation.)

  \begin{figure}
    \begin{center}
      \mettifig{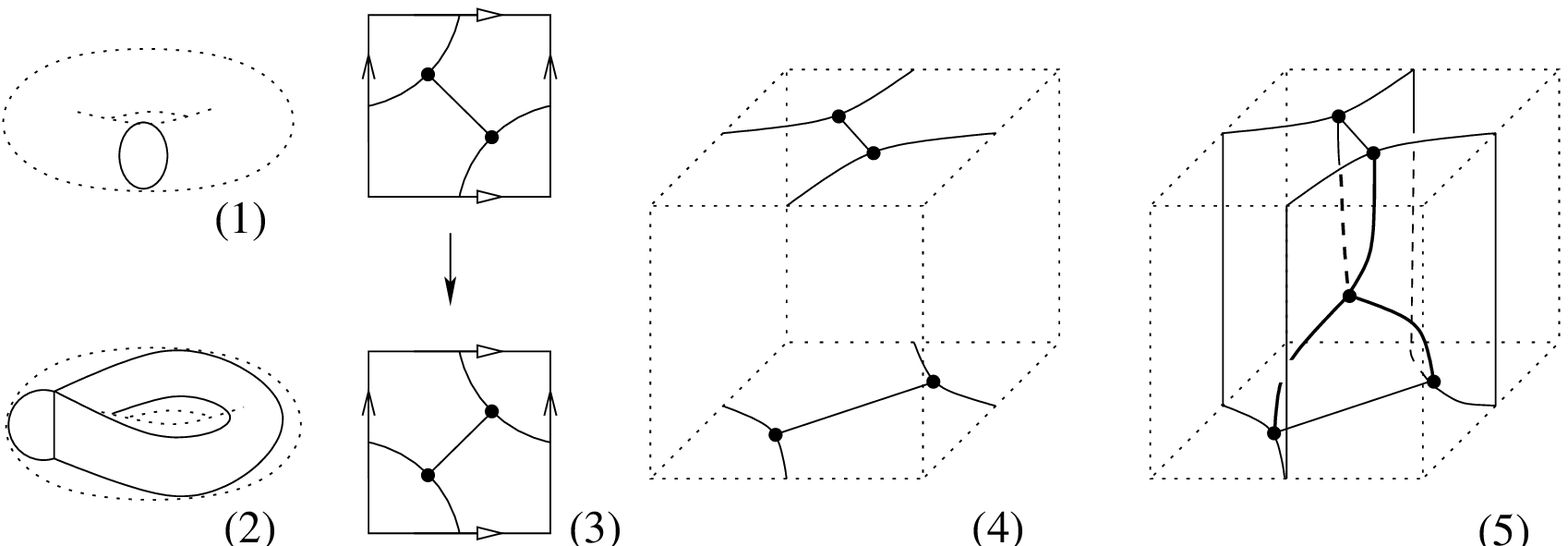, width = 14cm}
      \nota{Two polyhedra in the solid torus (1,2) and a flip (3,4) realized by adding a vertex (5)}
      \label{mobius_triplet:fig}
    \end{center}
  \end{figure}

The same technique, applied to both $F_1$ and $F_2$, allows to get a spine
of any $N(\alpha,\beta)$, as we will do below for $N(1,1)$, while
for $N(\alpha,\beta,\gamma)$ we will use a slightly different method.

\paragraph{Recognition of degenerate fillings}
  As we said, if $\alpha\in\{-2,-1,\infty\}$ we are in a special case, because $Q(\alpha)$
  does not contain the whole $\theta$-graph $\partial F_1$, but only a
  circle in it. The edge in $\partial F_1$ not contained in
  $Q(\alpha)$ is therefore adjacent to a single face $\Sigma$ in $P(\alpha)$,
  namely the face ($Q_1$, $E$ or $C$) of $P$ associated to $\alpha$.
We can therefore collapse $\Sigma$ and then further
collapse as much as possible, getting a polyhedron $P'(\alpha)$ which is still
  a spine of $N(\alpha)$. Along the collapse we kill some vertices of $P$, so
  $P'(\alpha)$ has less than $6$ vertices and can be easily analyzed.

To construct $P'(\alpha)$ we must first replace by three parallel strands the portion
of Fig.~\ref{P:fig}-(1) within a dotted rectangle and
then remove the attaching circle of the face $\Sigma$.
If $\Sigma$ is doubly incident
  to some edge $e$ then, after collapsing $\Sigma$,
    $e$ is again incident to one face $\Sigma'$ only, so we can collapse $\Sigma'$ too, and so on.
  The spine $P'(\alpha)$ we get here at the end of the process has a
    (maybe empty) singular set $S(P'(\alpha))$ which again is $4$-valent, but
  $P'(\alpha)\setminus S(P'(\alpha))$ may not consist of discs only,
    because some of the original discs of $P$ have been glued together, maybe
  producing annuli, M\"obius strips, or even more complicated surfaces.

  \begin{figure}
    \begin{center}
      \mettifig{Pinfty, width = 13cm}
      \nota{The spine $P'(\infty)$ is a torus}
      \label{Pinfty:fig}
    \end{center}
  \end{figure}

  \begin{figure}
    \begin{center}
      \mettifig{Pm1, width = 11cm}
      \nota{The spine $P'(-1)$, its modification, and the manifold it gives}
      \label{Pm1:fig}
    \end{center}
  \end{figure}

  For $\alpha = \infty$, everything collapses to the torus
  $T_2=\overline{F_2}$, as
  shown in Fig.~\ref{Pinfty:fig}, so $N(\infty) = T_2\times I$.
  For $\alpha=-1$, it is shown in
Fig.~\ref{Pm1:fig}-(1, 2) how $P'(-1)$, which has two vertices, can be transformed via an inverse
  ``lune move'' into a spine without vertices of the same manifold. This spine
  consists of a torus with a M\"obius strip attached, and it
  is a spine of $\seifuno A21$, since it is the inverse image in
  the Seifert fibration of $\seifuno A21$ of the graph \grafetto
  contained in the base annulus, as shown in
  Fig.~\ref{Pm1:fig}-(3). Note that the endpoint of \grafetto corresponds to the singular fibre.

  \begin{figure}
    \begin{center}
      \mettifig{Pm2, width = 11cm}
      \nota{The spine $P'(-2)$, and the manifold it gives}
      \label{Pm2:fig}
    \end{center}
  \end{figure}
  The spine $P'(-2)$ is shown in Fig.~\ref{Pm2:fig}-(1,2). It is made of one
  circle, one torus, and one annulus having
  one boundary component on the torus, and the other one turning three
  times along the circle. This is a spine of
  $\seifuno A31$ since it is the inverse image of the graph \grafetto
  in the base annulus of the fibration as in Fig.~\ref{Pm2:fig}-(3).

\paragraph{Next-to-degenerate fillings}
If $\alpha\in\{-3,-3/2,0\}$, the polyhedron $P'(\alpha)$ is
  constructed by attaching the M\"obius triplet to $P\setminus F_1$ along
  their $\theta$-shaped boundaries. Drawing $P'(\alpha)$ is easy:
  just replace the portion
of Fig.~\ref{P:fig}-(1) within a dotted rectangle by three ``monotonic'' strands.
Two of these strands must cross
  each other, while the other one must be straight vertical: the
  straight one is that representing the face associated
    to the slope in $\{-2,-1,\infty\}$ which is opposite to $\alpha$
    (\emph{i.e.}~the only slope having distance $2$ from $\alpha$).

  \begin{figure}
    \begin{center}
      \mettifig{Pm3, width = 11cm}
      \nota{The spine $P'(-3)$, its modification, and the manifold it gives}
      \label{Pm3:fig}
    \end{center}
  \end{figure}

The polyhedron
    $P'(-3)$ is as shown in Fig.~\ref{Pm3:fig}-(1), because $-3$ is opposite to
  $-1$, which corresponds to $Q_1$. Since $P'(-3)$ has one face with two vertices, we
  can apply to $P'(-3)$ a negative lune move,
getting a polyhedron which consists of a torus, a disc, and a
  M\"obius strip, as shown in Fig.~\ref{Pm3:fig}-(2).
    Actually, if we remove the M\"obius strip we get a
  torus and another M\"obius strip, as for $\alpha=-1$, whence a spine of
  $\seifuno A21$. The boundary loop of the removed M\"obius strip
  intersects the fibre of $\seifuno A21$ in one point, so it contributes with a
  $\seifuno D21$, with a vertical strip of $\seifuno D21$
  glued to a horizontal
  strip of $\seifuno A21$, as suggested by Fig.~\ref{Pm3:fig}-(3).
    The result is $N(-3)=\seifuno A21\bigu 0110\seifuno A21$.

\newcommand{\piccolaH}{
\begin{picture}(16,8)
\put(4,1){\line(0,1){6}}
\put(4,4){\line(1,0){8}}
\put(12,0){\line(0,1){8}}
\end{picture}
}

\newcommand{\staccalH}{
\begin{picture}(16,8)
\put(4,1){\line(0,1){2}}
\put(4,5){\line(0,1){2}}
\put(4,3){\line(1,0){8}}
\put(4,5){\line(1,0){8}}
\put(12,0){\line(0,1){8}}
\end{picture}
}

  \begin{figure}
    \begin{center}
      \mettifig{P0, width = 11cm}
      \nota{The spine $P'(0)$, its modification, and the manifold it gives}
      \label{P0:fig}
    \end{center}
  \end{figure}

The polyhedron $P'(0)$ is shown in Fig.~\ref{P0:fig}-(1) and
also in the top of part (2) of the same figure,
where each boundary arc should be glued
to some other one. We do not show
all the gluings, but we note that the left $\piccolaH$ is identified
``straight across'' with the $\piccolaH$ on the right.
The union of the two squares $U$ and $V$ gives an
annulus in $P'(0)$. We can now modify $P'(0)$
as suggested in Fig.~\ref{P0:fig}-(2), \emph{i.e.} we imagine
$\overline{U\cup V}$ to be made of two annular films glued together,
and we blow air between the films from the front of the figure.
The effect of this move on left and right boundary
graphs is given by \piccolaH$\to$\staccalH.
Using the exact identifications one sees that the new spine $R$, which
has two vertices, consists of one annulus and
  two discs. By removing the annulus, we get a polyhedron $Q$ without
  vertices, with $S(Q)=S^1$, and $Q\setminus S(Q)$
  consisting of one M\"obius strip whose boundary turns $3$ times along
  the circle $S(Q)$. This is a spine of $\seifdue D2131$, since it is the
  inverse image of a segment connecting the two singular points in the
  base disc, as suggested by Fig.~\ref{P0:fig}-(3).
The annulus $R\setminus Q$ we have removed will contribute with a $D\times S^1$
with two ``vertical'' (\emph{i.e.}~fibrewise) strips on
$\partial (D\times S^1)$ glued to $\partial \seifdue D2131$.
To analyze the gluing we must find which curves on
$\partial \seifdue D2131$ project in $Q$ to the boundary curves
$\gamma_0,\gamma_1$ of the annulus $R\setminus Q$.
Now $\gamma_0$ and $\gamma_1$ appear in the unfolded version
of $Q$ as in Fig.~\ref{P0:fig}-(3), and it is easy to see that
a curve on $\partial \seifdue D2131$ projecting in $Q$ to $\gamma_i$
must intersect the fibre once, so it has slope $1/k$ for some $k\in\matZ$
with respect to the natural homology basis. We deduce that
$$N(0)=\big(P\times S^1\big)\bigu 01{-1}k\seifdue D2131 .$$
Note that we cannot modify $k$ by
changing bases, because there is no free boundary component on the right block.
A more detailed analysis of $R$ would allow to compute the exact value of $k$,
but we will compute it later by an easier indirect argument.

\paragraph{Double filling}
Note that so far we have not used the flip of Fig.~\ref{mobius_triplet:fig}-(3,4,5),
but we do use it now to get a spine $P(1,1)$ for $N(1,1)$.
A neighbourhood of the singular set of $P(1,1)$ is obtained from
Fig.~\ref{P:fig}-(1) by removing both the portions within the dotted
boxes and replacing them by a ``curl,'' as shown in Fig.~\ref{P11:fig}-(1).
  \begin{figure}
    \begin{center}
      \mettifig{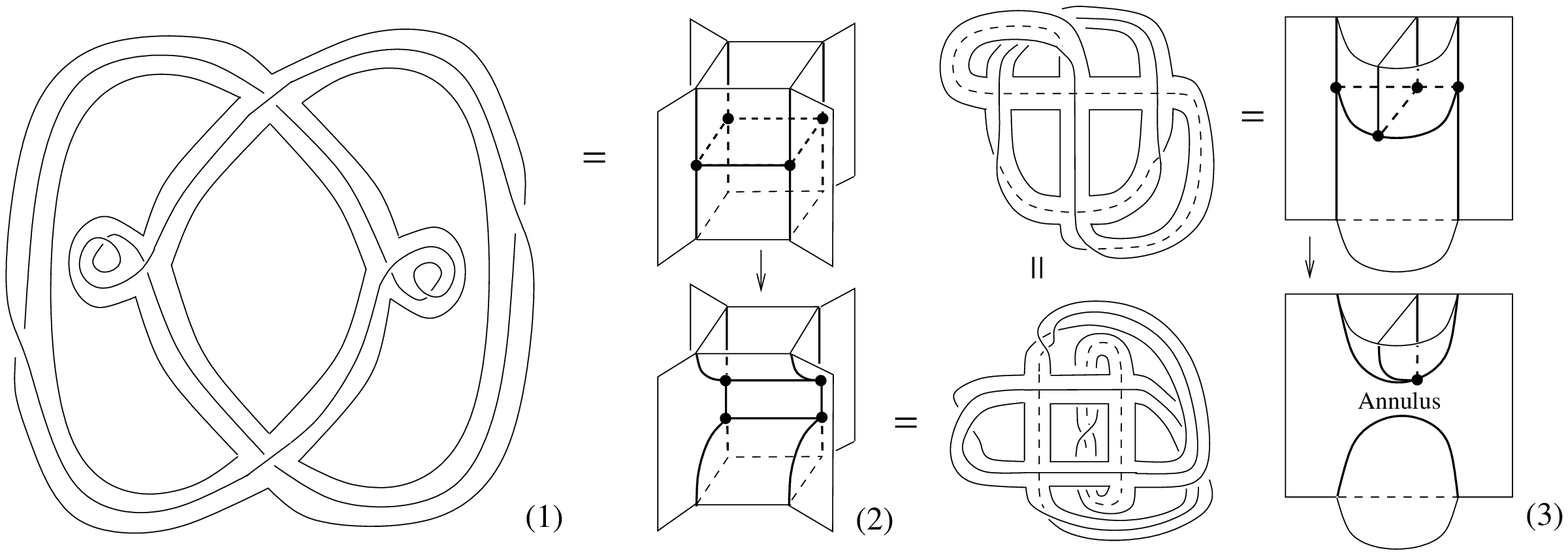, width = 11cm}
      \nota{The spine $P(1,1)$ and its modification}
      \label{P11:fig}
    \end{center}
  \end{figure}
The curl arises because the flip adds a vertex to the spine.
Parts (2) and (3) of Fig.~\ref{P11:fig} show how to transform $P(1,1)$ into
another spine $Q$ of $N(1,1)$ having only one vertex and
two faces (a disc and an annulus). The closure of
    the disc is a surface, which is easily seen to be a torus. Therefore $Q$
    is a torus with an annulus glued to it, and the two boundary circles
    of the annulus intersect in one point on the torus (giving the
    single vertex of $Q$). This implies that $N(1,1)=\big(P\times S^1\big)\bigb 0110$.

\paragraph{Recognition of $N(1,3,3)$ and $N(2,2,2)$}
To get a spine for $N(\alpha,\beta,\gamma)$, the above method
based on the spine $P$ of $N$ does not work.  We have used instead a
polyhedron $\widetilde{P}\subset N$ which is a spine of $N$ minus a ball
(not of $N$), and has $3$ hexagonal faces
$\widetilde{F}_1,\widetilde{F}_2,\widetilde{F}_3$
appearing in $\widetilde{P}$ as $F_1$ and $F_2$ appear in $P$.
Now the construction leading from $P$ to $P(\alpha)$ can be repeated for
all the $\widetilde{F}_i$'s and gives a spine
$\widetilde{P}(\alpha,\beta,\gamma)$ of
$N(\alpha,\beta,\gamma)$. The polyhedron $\widetilde{P}$ is more
complicated than $P$, and it is shown in~\cite{MaPe:1}
as a ``skeleton'' of the ``brick'' $B_9$.
Since $\widetilde{P}$ has $9$ vertices, $\widetilde{P}(1,3,3)$ has
$9+1+3+3 = 16$, while $\widetilde{P}(2,2,2)$ has $9+2+2+2=13$.
Recognizing $N(1,3,3)$ and $N(2,2,2)$ has therefore required several
steps, that we have decided to omit here.

\paragraph{Naturality of bases}
To show that the homology bases used so far in the proof are (up to sign) the meridian-longitude
bases, we state the following:

\begin{lemma}\label{exceptional:lemma}
With respect to the natural homology bases on $\partial N$,
the slopes $\alpha$ such that $N(\alpha)$ is non-hyperbolic
are $\infty,-3,-2,-1,0$.
\end{lemma}

This result can (and will) be proved very easily, using SnapPea for a finite
enumeration and the Gromov-Thurston $2\pi$-theorem.
Since these tools are part of a more elaborate argument that we use to conclude
the proof of Theorems~\ref{one:cusp:teo} to~\ref{three:cusps:teo}, we defer the
proof of Lemma~\ref{exceptional:lemma} to that argument, presented at the end of the paper.

Back to naturality, let
$h:\matQ\cup\{\infty\}\to\matQ\cup\{\infty\}$ be the map which
describes how the slopes change under the change of basis.
Lemma~\ref{exceptional:lemma} implies that $h$ sends the
set $\{\infty,-3,-2,-1,0\}$ to itself. Therefore it fixes $\infty$,
and it maps $t/u$ either to $t/u$ or to $-3-t/u$.
The above recognition of $N(0)$ and $N(-3)$, with coefficients
with respect to the basis coming from the spine $P$,
easily implies that $N(0)$
contains an incompressible annulus connecting two distinct boundary
components, while $N(-3)$ does not.
It is now easy to see that, with respect to the meridian-longitude basis, $N(0)$
does contain such an annulus: a spanning disc
of a component $C$ of the chain-link of Fig.~\ref{chain:fig}
gives a pair of pants $P$ in $N$ and, if we $0$-fill $C$,
the union of $P$ and a filling disc gives the desired annulus.
Therefore $h(0)=0$, and $h$ is the
identity, as required.

\paragraph{Conclusion for the slopes $-$1 and $-$2}
    We will now deduce relations (\ref{522:eqn}) to (\ref{526:eqn}) from
    (\ref{521:eqn}) and the topological informations on the
    fillings of $N$ just discovered.
    We start by proving (\ref{524:eqn}).
    Since $N(-1)=\seifuno A21$, using the fact that the boundary components
    of $N(-1)$  are symmetric, we deduce that there exists $\tiny\matr mnpq\in\GL_2(\matZ)$ such that
$$N\left(-1,\frac rs,\frac tu\right) = \seiftre {S^2}21{mr+ns}{pr+qs}{mt+nu}{pt+qu}.$$
Without loss of generality we can also assume $m\geqslant 0$. Since
$N(-1,\infty,\infty)$ equals $\seifdue {S^2}1011=S^3$,
we have $\seiftre {S^2}21mpmp = S^3$, which implies that $m=1$ and $p=0$.

Now recall that, among Seifert fibrations over $S^2$ with at most
$3$ singular fibres, only $3$ contain
non-separating tori, namely those described in
(\ref{from:Seif:to:T:eqn1}), (\ref{from:Seif:to:T:eqn2}), and (\ref{from:Seif:to:T:eqn3}).
We have shown that $N(1,1)$ contains a non-separating torus, so
$N(-1,1,1) = \seiftre{S^2}21{1+n}{q}{1+n}{q}$
also does, therefore it must be (up to orientation) one of the $3$ manifolds
just mentioned. We easily deduce that either $1+n = 4$ and $q=-1$,
whence relation (\ref{524:eqn}), or $1+n=-4$ and $q=1$.
The latter possibility is however absurd, because it would
imply that $N(-1,-2)$ is $\seifdue D217{-1}$, but we know that
$N(-2)=\seifuno A31$.

Relation (\ref{523:eqn}) is
    proved in the same way. We first note that
  $$N\left(-2,\frac rs,\frac tu\right) = \seiftre {S^2}32{mr+ns}{pr+qs}{mt+nu}{pt+qu}.$$
Using the identity
$N(-2,\infty,\infty)=S^3$ we then get $m=1$ and $p=0$. Moreover
    $N(-2,1,1)$ contains a non-separating torus, so $1+n = 3$
    and $q=-1$.

\paragraph{Conclusion for the slope $-$3}
    Concerning (\ref{522:eqn}), we have
  $$N\left(-3,\frac rs,\frac tu\right) =
    \seifdue D21{mr+ns}{pr+qs}\bigu 0110 \seifdue D21{mt+nu}{pt+qu}$$
    for some $\tiny\matr mnpq\in\GL_2(\matZ)$ with $m\geqslant 0$.
    We first note that by reversing orientations we get two blocks of
    type
    $$\seifdue D2{-1}{mr+ns}{-pr-qs} = \seifdue
    D21{mr+ns}{-(m+p)r-(n+q)s},$$
    glued along the same $\tiny\matr 0110$. Since $
    \det\tiny\matr mn{-(m+p)}{-(n+q)}=-\det\tiny\matr mnpq$, we can
    suppose $\det\tiny\matr mnpq = 1$.
    Now, we know that
$$\lens 31=N(-3,\infty,\infty)=\seifdue D21mp\bigu 0110 \seifdue D21mp.$$
An obviously necessary condition for this to hold is that $m=1$.
In addition
$\seifuno D2{1+2p}\bigu 0110\seifuno D2{1+2p} = \seifdue
{S^2}2{1+2p}{-1-2p}2=\lens{4-(1+2p)^2}{p(1+2p)-2}$
must be equal to $L(3,1)$, whence $4-(1+2p)^2=\pm 3$.
This implies that $p\in\{0,-1\}$.

Now we use the fact that
$N(-3,-2)$ and $N(-3,-1)$, being single fillings of $\seifuno A21$ and
$\seifuno A31$ respectively, are either reducible or atoroidal.
However
\begin{eqnarray*}
N(-3,-2) &=& \seifdue D21{2-n}{2p-q} \bigu 0110 \seifuno A21\\
N(-3,-1) &=& \seifdue D21{1-n}{p-q} \bigu 0110 \seifuno A21
\end{eqnarray*}
so in both these expressions the block on the right-hand side must be
either reducible or a solid torus. This implies that $|2-n|\leqslant 1$ and $|1-n|\leqslant 1$,
whence $n\in\{1,2\}$.

To conclude we examine $N(-3,1,1)$ and note that it can be expressed both
as $\seifuno Aab\bigb 0110$,
for suitable $a$ and $b$, and as
$$\seifdue D21{1+n}{p+q} \bigu 0110 \seifdue D21{1+n}{p+q}.$$
Uniqueness of the JSJ decomposition implies that one block of this
expression, and hence also the other one, must have a non-unique fibration.
So $n=1$, in which case $\seifdue D2121$ can be fibred also as
$S\timtil S^1$. In addition, the alternative fibration
of one block must be matched to either the original or the alternative
fibration of the other block. Using relations~(\ref{from:D2X:to:D2X:eqn}),
(\ref{from:D2X:to:SX:eqn}), and (\ref{from:XD2:to:XS:eqn}) one sees that
this can only happen for $(p,q)=(0,1)$ and $(p,q)=(-1,-2)$, but
${\rm det} {\tiny\matr 11pq}=1$, so $(p,q)=(0,1)$, as
required. Even if not necessary, note that we have shown the equality
$N(-3,1,1)=(K,1)$, which is coherent with~(\ref{526:eqn})
and~(\ref{from:Seif:to:T:eqn4}).

\paragraph{Conclusion for the slope 0}
Concerning (\ref{525:eqn}), we have
$$N\left(0,\frac rs,\frac tu\right) = \seifdue
D{mr+ns}{pr+qs}{mt+nu}{pt+qu} \bigu 01{-1}k \seifdue D2131$$
for some $\tiny\matr  mnpq\in\GL_2(\matZ)$ with
$p\geqslant 0$ and some $k\in\matZ$.
Since $S^2\times S^1 = N(0,\infty,\infty) = \seifdue Dmpmp \bigu 01{-1}k \seifdue D2131$,
we have $m=0$, $p=1$, and
$$S^2\times S^1 =
(S^2\times S^1)\#\seiftre {S^2}21311k.$$
Therefore
$\seifdue {S^2}2{1+2k}31 = S^3$, which implies
that $5+6k=\pm 1$, hence $k=-1$.

Having shown that $p=1$, we have $n=\pm 1$.
We prove by contradiction that $n=1$. Assuming $n=-1$ we get
    \begin{eqnarray*}
    N(0,-2,t/u) & = & \seifdue D1{2-q}{-u}{t+qu}\bigu 01{-1}{-1} \seifdue D2131 \\
    & = & \seiftre {S^2}2{-1}31{t+2(q-1)u}{-u}.
    \end{eqnarray*}
Equation~(\ref{523:eqn}), now proved, shows that for all $t/u$
this manifold must be equal to $\seiftre {S^2}2{-1}32{t+2u}{-u}$,
which is clearly impossible whatever $q$.

Having shown that $n=1$, we determine $q$ noting that
    \begin{eqnarray*}
    N(0,1,1) & = & \seifdue D1{1+q}1{1+q} \bigu 01{-1}{-1} \seifdue D2131 \\
    & = & \seiftreeul {S^2}2131{2(1+q)}1{-1}
    \end{eqnarray*}
contains a non-separating torus, hence $q=2$, as required.

\paragraph{Conclusion for the slopes (1,1)}
Finally, we prove~(\ref{526:eqn}). We have
$$N\left(1,1,\frac tu\right) = \seifuno A{mt+nu}{pt+qu}/_{\tiny\matr 0110}$$
for some $\tiny\matr  mnpq\in\GL_2(\matZ)$ with $p\geqslant 0$.
First, $S^2\times S^1 = N(1,1,\infty)$
$= \seifuno Amp/_{\tiny\matr 0110}$ gives $m=0$ and $p=1$, hence
$n=\pm 1$. We also deduce from relations~(\ref{522:eqn})
to~(\ref{525:eqn})
and~(\ref{from:Seif:to:T:eqn1}) to~(\ref{from:Seif:to:T:eqn4})
that $N(1,1,\alpha) = T_{\tiny\matr {\alpha+1} 1{-1}0}$
for $\alpha\in\{-3,-2,-1,0\}$.
Using~(\ref{from:A:to:T:eqn}) we easily
get $q=1$ and $n=1$, as required.
\end{proof}

We can now prove our main results.

\vspace{2pt}\noindent\textit{Proof of \ref{one:cusp:teo},
\ref{two:cusps:teo}, and \ref{three:cusps:teo}\ }
The proof is divided in $4$ steps:
\begin{enumerate}
\item
we use
Proposition~\ref{recognized:fillings:prop} and Lemma~\ref{graph:lem}
to recognize the JSJ decomposition of many non-hyperbolic fillings;
\item
we apply the symmetries listed in
Proposition~\ref{repetitions:prop} finding more
non-hyperbolic fillings;
\item
we use the Gromov-Thurston $2\pi$-theorem
to prove that $M$ is negatively curved except in the cases listed
(up to permutation) in
Steps 1 and 2 and in a finite number of other cases.
As already mentioned, the orbifold
theorem shows that $M$ is hyperbolic if it is negatively curved;
\item
using SnapPea, we examine the finite number of cases
not covered by the previous steps,
and we check that the corresponding $M$ is always hyperbolic.
\end{enumerate}

\paragraph{Step 1: Recognition of non-hyperbolic fillings}
The equalities listed in Theorem~\ref{one:cusp:teo} are deduced from
Proposition~\ref{recognized:fillings:prop}. We now prove most of the
equalities listed in Tables~\ref{two:cusps:table},~\ref{closed:fill:table1},
and~\ref{closed:fill:table2}, by studying the ``degenerate'' cases in which the
expressions in Proposition~\ref{recognized:fillings:prop} do not
describe the JSJ decomposition of the manifold.
We say that $\big(F,(\scriptstyle{i_1},\scriptstyle{j_1}),\ldots,
(\scriptstyle{i_k},\scriptstyle{j_k})\big)$ is \emph{genuine} if
it is irreducible and not a solid torus. Degenerate cases
occur mostly in presence of non-genuine blocks.

\paragraph{1-cusped fillings}
We first prove all the equalities listed in Table~\ref{two:cusps:table}, except
the last two. For $(p/q,r/s)=(1,1)$ the equality follows from
Proposition~\ref{recognized:fillings:prop} directly.
For $p/q\in\{-3,-2,-1,0\}$ the argument is quite simple: we only
spell it out here for $p/q=-3$, the other cases being even easier.

So, assume $p/q=-3$ and $r/s\in\matQ$. By~(\ref{522:eqn}) we have
$$M=\seifdue D21{r+s}s\bigu0110\seifuno A21.$$
this formula describes the JSJ decomposition of $M$ except
in the following cases:
\begin{description}
\item[(A)] the block on the left-hand side
is non-genuine;
\item[(B)] the block on the left-hand side
admits an alternative fibration which is matched to the fibration
of the block on the right-hand side.
\end{description}
Case (A) occurs for $|r+s|\leqslant 1$.
When $|r+s|=0$, \emph{i.e.}~$r/s=-1$, we have $M=\matRP^3\#\seifuno
D21=\matRP^3\#(D\times S^1)$ by (\ref{from:DX:to:lensX:eqn}).
If $|r+s|=1$, \emph{i.e.}~$r/s=-1+\frac1n$ for some
$n\in\matZ\setminus\{0\}$, we have
$$M=\seifuno D1{1+2n}\bigu0110 \seifuno A21 = \seifdue D{1+2n}{-2}21 =
\seifdue D{1+2n}221$$
by (\ref{from:D1D2:to:D2:eqn}), (\ref{minus:eqn}), and
(\ref{from:F:to:F:eqn}), which is genuine except for $|1+2n|\leqslant
1$, which only occurs for $n=-1$, giving $D\times S^1$.

Using Lemma~\ref{graph:lem} it is not hard to see that case (B)
occurs only for $r/s=-3$, in which case $M=\seifuno S21$.

We now prove all the equalities in Table~\ref{closed:fill:table1}, and the
equalities in Table~\ref{closed:fill:table2} corresponding to the case $(p/q,r/s)=(1,1)$
\paragraph{The slope $-$3} We start by considering the case where $p/q=-3$.
Relation~(\ref{522:eqn}) then shows that $M$ is homeomorphic to
\begin{equation}\label{522:at:work:eqn}
\seifdue D21{r+s}s\bigu 0110 \seifdue D21{t+u}u
\end{equation}
The first case where the block on the left-hand side
is non-genuine occurs when $r/s=-1$, in which case the block becomes
$\seifdue D2101$ and we deduce from (\ref{from:DX:to:lensX:eqn})
and (\ref{from:D2:to:D1:eqn}) that
$$M=\matRP^3\#\seifdue{S^2}21{t+u}u.$$
Now we employ~(\ref{from:D1:to:lens:eqn})
to see that the second summand is $\lens{t+3u}u$.
Summing up, we have  $M=\matRP^3\#\lens{t+3u}u$ when $r/s=-1$ (and $p/q=-3$).

Another case where the block on left-hand side of~(\ref{522:at:work:eqn})
is non-genuine occurs when $|r+s|=1$, namely when $r/s=-1+1/n$ for
some $n\in\matZ$, $n\ne 0$. In this case the block is
$\seifdue D211n = \seifuno D2{1+2n}$. We then use~(\ref{from:D1D2:to:D2:eqn})
to deduce that
$$M=\seiftre{S^2}{2n+1}{-2}21{t+u}u$$
in this case. Now this manifold is Seifert with 3 exceptional fibres unless
$|2n+1|\leqslant 1$ or $|t+u|\leqslant 1$. The former inequality is only possible
for $n=-1$ (recall that $n\ne 0$), \emph{i.e.}~for $r/s=-2$, and
$M=\seiftre {S^2}21{t+u}u12 =\seifdue {S^2}25{t+u}u$ in this case. Now we
apply~(\ref{from:D1:to:lens:eqn}) to see that
$$M=\lens{2u+5(t+u)}{u+2(t+u)}=\lens{5t+7u}{2t+3u}.$$
Let us assume instead that $|t+u|\leqslant 1$. We dismiss the case $t/u=-1$
which was already considered up to permutation, so we have
$t/u=-1+1/m$ for $m\in\matZ$, $m\ne 0$, and
$M=\seiftre {S^2}211m{2m+1}{-2} = \seifdue {S^2}2{2m+1}{2m+1}{-2}$.
We now apply~(\ref{from:D1:to:lens:eqn}), eventually getting
$$M=\lens{(2n+1)(2m+1)-4}{(2n+1)m-2}.$$
We also note that the lack of symmetry in $n$ and $m$ is only apparent,
because
$$(2nm+m-2)(2nm+n-2)\equiv 1 \quad({\rm mod}\ (2n+1)(2m+1)-4).$$
We have now recognized $M$ whenever $r/s$ equals $-1$ or $-1+1/n$
(and $p/q=-3$). So we assume neither $r/s$ nor $t/u$ have this
form, \emph{i.e.}~$|r+s|>1$ and $|t+u|>1$. So~(\ref{522:at:work:eqn})
describes the gluing of two genuine Seifert blocks, and the gluing
matrix does not match the fibres of the fibrations implicit in the
expressions of the blocks. Therefore~(\ref{522:at:work:eqn})
gives the JSJ decomposition of $M$ except in the following cases (up to
permutation):
\begin{itemize}
\item[(A)] the block on left-hand side has an alternative fibration
which is matched to the original fibration of the block on the right-hand side;
\item[(B)] both the blocks admit alternative fibrations, and these
fibrations are matched.
\end{itemize}
Now the left-hand side has an alternative fibration precisely when
$|r+s|=2$, namely when $r/s=-1+\frac{2}{2n+1}$ for some $n\in\matZ$,
and in this case we deduce from (\ref{from:D2X:to:D2X:eqn}) and~(\ref{from:D2X:to:SX:eqn}) that
\begin{eqnarray*}
M & = & \seifdue D212{2n+1}\bigu 0110 \seifdue D21{t+u}u \\
& = & \seifdue D2121\bigu {-n}110 \seifdue D21{t+u}u \\
& = & \big(S\timtil S^1\big)\bigu 1{n+1}0{-1} \seifdue D21{t+u}u.
\end{eqnarray*}
Now case (A) occurs precisely when $n+1=0$, whence $r/s=-3$, and
$$M = \big(S\timtil S^1\big)\bigu 100{-1} \seifdue D21{t+u}u.$$
Up to changing sign the gluing matrix is now properly matching both the sections
and the fibres, so $M$ is Seifert fibred over the gluing of $S$ and $D$, which
gives
$$M=\seifdue {\matRP^2}21{t+u}u$$
for $r/s=-3$ (and $p/q=-3$).

Assume case (B) occurs. Then $t/u$ also has the form $-1+\frac{2}{2m+1}$,
and by (\ref{from:XD2:to:XD2:eqn}) and~(\ref{from:XD2:to:XS:eqn}) we have
\begin{eqnarray*}
M & = & \big(S\timtil S^1\big)\bigu 1{n+1}0{-1} \seifdue D212{2m+1} \\
& = & \big(S\timtil S^1\big)\bigu 1{n+1}m{mn+m-1} \seifdue D2121 \\
& = & \big(S\timtil S^1\big)\bigu {m+1}{m+n+mn}{-1}{-n-1}
\big(S\timtil S^1\big).
\end{eqnarray*}
Moreover the fibres should match, so $m+n+mn=0$. Therefore we have either
$n=m=0$  or $n=m=-2$, namely either $r/s=t/u=1$ or
$r/s=t/u=-5/3$. The two manifolds we get are
$$M=\big(S\timtil S^1\big)
\bigcup\nolimits_{\tiny{\matr {-1}0{\pm 1}1}}
\big(S\timtil S^1\big)$$
which are easily seen to be $(K,\pm 1)$, where $K$ is the Klein bottle.
Finally, we have $(K,+1)=(K,-1)$.

\paragraph{The slope $-$2} Let us consider now the case $p/q=-2$,
where equation~(\ref{523:eqn}) shows that $M$ is
$$\seiftre {S^2}32{r+2s}{-s}{t+2u}{-u}.$$
There are three exceptional fibres unless $|r+2s|\leqslant  1$
(up to permutation).
If $r/s=-2$, we have $M=\seiftre{S^2}3201{t+2u}{-u}$, which is equal
to $\lens 31\#\lens {t+2u}u$ by~(\ref{from:D2:to:lens:plus:lens:eqn}).
If $|r+2s|=1$, \emph{i.e.}~$r/s=-2+1/n$, then
$M = \seiftre{S^2} 321{-n}{t+2u}{-u} = \seifdue {S^2}3{2-3n}{t+2u}{-u}$. We
apply~(\ref{from:D1:to:lens:eqn}), getting
$$M=\lens{3u-(2-3n)(t+2u)}{u-(1-n)(t+2u)}=\lens{3n(t+2u)-2t-u}{n(t+2u)-t-u}.$$

\paragraph{The slope $-$1} The case $p/q=-1$ is similar to the
previous one. Relation~(\ref{524:eqn}) proves that $M$ is
$$\seiftre {S^2}21{r+3s}{-s}{t+3u}{-u}.$$
This manifold is Seifert fibred over $S^2$, and there are three exceptional fibres
unless $|r+3s|\leqslant  1$ up to permutation.
If $r/s=-3$ then $M=\seiftre{S^2}2110{t+3u}{-u}$, which is equal to
$\matRP^3\#\lens{t+3u}u$ by~(\ref{from:D2:to:lens:plus:lens:eqn}).
If $|r+3s|=1$, \emph{i.e.}~$r/s=-3+1/n$, we have
$M=\seiftre{S^2}211{-n}{t+3u}{-u} = \seifdue{S^2}2{1-2n}{t+3u}{-u}$.
Now we use~(\ref{from:D1:to:lens:eqn}), deducing that
$$M=\lens{-2u+(1-2n)(t+3u)}{u+n(t+3u)} = \lens{2n(t+3u)-t-u}{n(t+3u)+u}.$$

\paragraph{The slope 0} Let us turn now to the case $p/q=0$, where~(\ref{525:eqn}) shows that
$M$ is
\begin{equation}\label{525:at:work:eqn}
\seifdue Ds{r+2s}u{t+2u}\bigu 0{-1}11 \seifdue D2131
\end{equation}
The right-hand side block is now always genuine, but the left-hand side block
is not if $|s|=1$ (recall that $s\ne 0$), namely when
$r/s=n\in\matZ$. In this case the block becomes
$\seifdue D1{n+2}u{t+2u} = \seifuno Du{t+(n+4)u}$.
Now we use~(\ref{from:D1D2:to:D2:eqn}) to see that
$$M=\seiftre{S^2}2131{t+(4+n)u}{-t-(3+n)u} = \seiftre{S^2}2{-1}31{t+(4+n)u}u.$$
So $M$ is Seifert with three exceptional fibres
unless $|t+(4+n)u|\leqslant  1$.
If $t+(4+n)u=0$, namely $t/u=-4-n$, we
use~(\ref{from:D2:to:lens:plus:lens:eqn}) and see that
$M=\matRP^3\# \lens 31$. If $|t+(4+n)u|=1$, \emph{i.e.}~$t/u=-4-n+1/m$,
we have $M=\seiftre {S^2}21311{m-1} = \seifdue {S^2}312{2m-1}$.
We eventually use~(\ref{from:D1:to:lens:eqn}) and see that $M=\lens{6m-1}{2m-1}$.

We consider now the case where
$r/s,t/u\not\in\matZ$ in~(\ref{525:at:work:eqn}), so both the blocks involved
are genuine Seifert, and the gluing matrix does not match the fibrations.
In addition, the right-hand side block has a unique fibration, so the JSJ
decomposition of $M$ is given by
\begin{equation}\label{525:JSJ:eqn}
\seifdue Ds{r+2s}u{t+2u}
\bigcup\nolimits_{\tiny{\matr 0{-1}11}} \seifdue D2131
\end{equation}
unless the left-hand side block has an alternative fibration
whose fibre is matched to the fibre of $(D,(2,1),(3,1))$.
Assuming this happens, we must have $|s|=|u|=2$, so
$r/s=n+1/2$ and $t/u=m+1/2$.
By (\ref{from:D2X:to:D2X:eqn}) and~(\ref{from:D2X:to:SX:eqn}) we then have
\begin{eqnarray*}
M & = & \seifdue D2{2n+5}2{2m+5} \bigu 0{-1}11 \seifdue D2131 \\
 & = & \seifdue D2121 \bigu {m+n+4}{-1}{-m-n-3}1 \seifdue D2131 \\
 & = & \big(S\timtil S^1\big) \bigu {-1}{-m-n-5}1{m+n+4} \seifdue D2131.
\end{eqnarray*}
Fibres match for $n+m+5=0$, \emph{i.e.}~$m=-n-5$.
Summing up, we have shown that for $r/s,t/u\not\in\matZ$ the JSJ
decomposition of $M$ is given by~(\ref{525:JSJ:eqn}) unless
$r/s=n+1/2$ and $t/u=-n-9/2$. In this case we have,
using~(\ref{from:XD2:to:XD2:eqn}) with $l=-1$ and switching the orientation of
$S\timtil S^1$:
\begin{eqnarray*}
M &=& \big(S\timtil S^1\big)\bigcup\nolimits_{\tiny{\matr {-1}01{-1}}}
  \seifdue D2131 \\
 &=&\big(S\timtil S^1\big)\bigcup\nolimits_{\tiny{\matr {-1}00{-1}}}
  \seifdue D2{-1}31 \\
 &=&\big(S\timtil S^1\big)\bigcup\nolimits_{\tiny{\matr {-1}001}}
    \seifdue D2{-1}31 \\
 &=&\seifdueeul {\matRP^2}2131{-1}.
\end{eqnarray*}

\paragraph{The slopes (1,1)}
We now study the case $p/q=r/s=1$, where equation~(\ref{526:eqn}) shows that
$$M = \seifuno Au{t+u}\bigb 0110.$$
The block which is being glued to itself is non-genuine if
$t/u=n\in\matZ$. If this is the case the block becomes
$\seifuno A1{n+1}$ and we deduce from~(\ref{from:D2X:to:D2X:eqn}) and~(\ref{from:A:to:T:eqn}) that
$$M  =  \big(A\times S^1\big)\Big/_{\tiny{\matr{-1-n}110}} \\
=  T_{\tiny{\matr{n+1}1{-1}0}}$$

If $t/u\not\in\matZ$ then $\seifuno A{u}{t+u}$ is genuine Seifert
with a unique fibration and
the self-gluing does not match the fibres, so the JSJ decomposition of $M$
is as required.

\paragraph{Step 2: Recognition of more non-hyperbolic fillings}
As announced, we use now relations~(\ref{511right:eqn})
to~(\ref{515:eqn}) to
find triples which were not listed in Step~1 but give the same fillings
as triples listed in Step~1. This will complete the proof of all the
equalities listed in Tables~\ref{two:cusps:table} to~\ref{closed:fill:table3}.
We can immediately conclude with Table~\ref{two:cusps:table}:
Proposition~\ref{repetitions:prop} gives
$N(-3/2,-5/2)=N(-3/2,0)=N(-4,-1/2)$,
therefore the last two lines of Table~\ref{two:cusps:table} are
deduced from the case $(p/q,r/s)=(0,-3/2)$ in the same table.

All equalities of table~\ref{closed:fill:table1} have already been proved.
We therefore turn to Tables~\ref{closed:fill:table2} and~\ref{closed:fill:table3}.
For the sake of conciseness and clarity,
we have organized our constructions in Table~\ref{table6}.
\begin{table}\begin{center}
\begin{tabular}{c|c|c}
Original triple & Move & New triple \\ \hline\hline
$-\frac 32,-\frac 52,\frac tu\phantom{\Big|}$   &
    (\ref{511left:eqn}) &
    $-\frac 32,0,-\frac u{t+2u}-2$  \\ \hline
$-4,-\frac 12,\frac tu\phantom{\Big|}$  &
    inverse of (\ref{511right:eqn}) &
    $-\frac 32,0,-\frac tu-3$   \\ \hline
$-\frac 32,-\frac 73,-\frac 73\phantom{\Big|}$  &
    (\ref{511left:eqn}) &
    $-\frac 32,1,1$     \\ \hline
$-4,-4,-\frac 23\phantom{\Big|}$    &
    inverse of (\ref{511right:eqn}) &
    $-\frac 32,1,1$     \\ \hline
$-\frac 32,-\frac 32,-4\phantom{\Big|}$     &
    (\ref{511right:eqn}) &
    $-4,1,1$    \\ \hline
$-\frac 52,-\frac 52,-\frac 43\phantom{\Big|}$  &
    (\ref{512:eqn}) &
    $-\frac 52,1,1$     \\ \hline
$-\frac 12,-5,-5\phantom{\Big|}$    &
    (\ref{513:eqn}) &
    $-\frac 12,1,1$     \\ \hline
$\begin{array}{c}1,2,n\\ n=2,3,4,5\end{array}$  &
    (\ref{514:eqn}) &
    $1,2,2-n$   \\ \hline
$1,-4,-\frac 13\phantom{\Big|}$     &
    (\ref{515:eqn}) &
    $1,-4,-3$   \\ \hline\hline
$-\frac 32,-\frac 32,-\frac 83\phantom{\Big|}$  &
    (\ref{511right:eqn}) &
    $-4,-\frac 13,1$    \\ \hline
$-\frac 52,-\frac 53,-\frac 53\phantom{\Big|}$  &
    (\ref{512:eqn}) &
    $-\frac 52,-4, -\frac 12$
\end{tabular}\end{center}
\nota{Sporadic triples giving non-hyperbolic fillings}\label{table6}
\end{table}
Let us first comment on the first nine lines of this table.
Here we find on the first column a triple not found in Step~1 (except
for $t/u\in\{-3,-2,-1,0\}$ in the first two lines).
In the second column there is
a reference to one of the symmetry relations listed in Proposition~\ref{repetitions:prop}.
We view this relation
as a move acting on the triple to its left, and in the third column we give
the triple resulting from this action; this triple
is readily seen to be one of those of Step~1.
The last two lines are similar, except that in the last column
the triple is \emph{not} one of those of Step~1. Instead, it
appears in the first column of the ninth and second line respectively.

Let us now recognize $M$ in all the lines listed in the first column of Table~\ref{table6}.
Of course this is trivial for all the lines except the first two,
because for the other cases we only
have to use Step~1 for the triple on the last column (but in Table~\ref{closed:fill:table3}
we have often simplified the expression using formulae~(\ref{from:D2:to:D1:eqn}),
(\ref{from:D2X:to:D2X:eqn}), and~(\ref{from:XD2:to:XD2:eqn}), and changes of orientation).

Then, let us start from the
case $-\frac 32,-\frac 52,\frac tu$, so $M$ is the same as the manifold obtained with the
filling coefficients $0,-\frac 32,-\frac {2t+5u}{t+2u}$.
First, if $t+2u=0$, \emph{i.e.}~$t/u=-2$, we easily get from the case
$p/q=\infty$ that $M=\lens 21=\matRP^3$. Otherwise we have to refer to the
case $p/q=0$. Now $-\frac {2t+5u}{t+2u}$ is an integer when
$|t+2u|=1$, \emph{i.e.}~$t/u=-2+1/n$, and in this case
$-\frac {2t+5u}{t+2u}=-n-2$.  We further need to make distinctions
depending on whether
$$-3+(4-(n+2))2=-2n+1$$
has value $0$ or $\pm1$. Of course $0$ cannot be attained, and $+1$
also cannot (because $n\ne 0$), so the only special case occurs for $n=1$,
namely for $t/u=-1$, and $M=\lens{13}5$, as we knew already from the discussion of slope $-1$.
We then deduce that for $t/u=-2+1/n$ but $n\ne 1$ we have
$$M=\seiftre {S^2}{2}{-1}{3}{1}{3-(4-n-2)2}{-2}
=\seiftre {S^2}{2}{-1}{3}{1}{2n-1}{-2}.$$
When $-\frac{2t+5u}{t+2u}$ is neither an integer nor infinity,
\emph{i.e.}~when $t/u$ is not $-2$ or $-2+1/n$, we will
have, after easy calculations,
$$M=\seifdue D2131\bigcup\nolimits_{\tiny{\matr 110{-1}}}
            \seifdue D21{t+2u}u\big)$$
except in the special case where
$-\frac {2t+5u}{t+2u}=-(-\frac 32)-4$, which occurs for $t/u=0$, and
$$M=\seifdueeul {\matRP^2}2131{-1}$$
in this case, as we knew already.

The discussion of the cases $-4,-1/2,t/u$ is easier, because
it reduces to that of $0,-3/2,-t/u-3$. Namely, for $t/u=n$ we have
to consider the special cases $n=0,-1$, leading to
$\lens{13}5$ and $\lens{11}3$ respectively (as we knew already), otherwise
$$M=\seiftre {S^2}2{-1}31{2n+1}{-2}.$$
Moreover we have the exceptional case where $-t/u-3=-(-3/2)-4$, which occurs
for $t/u=-1/2$ and leads to
$$M=\seifdueeul {\matRP^2}2131{-1}.$$
In all other cases we have
$$M=\seifdue D2131\bigu 1110 \seifdue D21ut.$$

Now that we have recognized the $M$'s arising from the triples on
the first column of Table~\ref{table6}, we could show, with a little patience,
that even a repeated application of the
symmetries~(\ref{511right:eqn}) to~(\ref{515:eqn})
to any of the triples
found so far does not lead to new triples.  This will actually
be deduced from the rest of our argument, because we will see that all the
triples, except those listed so far (up to permutation),
give a hyperbolic $M$.

\paragraph{Step 3: Application of the $2\pi$-theorem}
As announced, we will now use the Gromov-Thurston $2\pi$-theorem~\cite{2pi:cite}
to show that  a Dehn filling of $N$ is negatively curved for all fillings
except those listed in Steps~1 and~2 and a finite number of other ones.
More precisely, let us define $\calG$ as the
set of all slopes, pairs of slopes, and triples
of slopes that were proved in Steps~1 and~2 to give a non hyperbolic
filling of $N$. The aim of Step~3 is to construct three finite sets of slopes
$S_1,S_2,S_3$ such that \emph{if $\alpha,\beta,\gamma$ are slopes on $\partial N$
and $N(\alpha,\beta,\gamma)$ is not negatively curved then either
$(\alpha,\beta,\gamma)\in\calG$ or, up to permutation,
$(\alpha,\beta,\gamma)\in S_1\times S_2\times S_3$.}

We begin by
recalling the statement of the $2\pi$-theorem, so we consider a cusped hyperbolic
manifold $M$ and a horospherical cusp section $H$ of $M$, that is an embedded
surface parallel to $\partial M$ such that the metric of $M$ restricts to a Euclidean
metric on $H$. Then $H=\sqcup_{i=1}^kT'_i$, where $T'_i$ is a Euclidean torus isotopic
to the $i$-th component $T_i$ of $\partial M$. If $\alpha_i$ is a slope
on $T_i$ we now define the length of $\alpha_i$ with respect to $H$
(or $H$-length for short) as the
shortest length of a curve in $T'_i$ parallel to $\alpha_i$. The Gromov-Thurston
theorem now asserts that $M(\alpha_{i_1},\ldots,\alpha_{i_h})$ is negatively
curved whenever all filling slopes $\alpha_{i_j}$ have $H$-length larger than $2\pi$.

Before proceeding we note that if a homology basis $(a,b)$ is fixed on a torus $T$,
the Euclidean structures on $T$ up to isotopy are determined by two parameters
$x+iy\in\matC$ with $y>0$ and $A\in\matR$ with $A>0$. Here $x+iy$ gives
the \emph{shape} of $T$ (\emph{i.e.}~the Euclidean structure up to scaling)
while $A$ gives the area. The relation with the fixed basis $(a,b)$ comes from the
fact that there should exist a locally isometric covering
$\matC\to T$ on which $a$ and $b$ act as the translations
$z\mapsto z+\sqrt{A/y}$ and $z\mapsto z+(x+iy)\sqrt{A/y}$ respectively.
It is now an easy exercise to show that the minimal length $\ell(p/q)$
of a slope $pa+qb$ on $T$ satisfies
\begin{equation}\label{ell:eqn}
\ell(p/q)^2=\frac Ay\big((p+xq)^2+(yq)^2\big).
\end{equation}

We now specialize to our chain-link exterior $N$ and note that it has symmetries
realizing any permutation of the boundary components, so the shapes
of its three cusps are actually the same. We describe this common shape using
the homology bases $(a_i,b_i)$ chosen by SnapPea, where
$a_i=\mu_i$ and $b_i=2\mu_i-\lambda_i$ for $i=1,2,3$, with
$(\mu_i,\lambda_i)$ being the meridian-longitude basis used above.
The shape is then described by the parameter
$1/2 +i\sqrt 7/2$.
Moreover, if we choose a maximal horospherical
cusp section $H$ consisting of three tori having the same area, this area turns
out to be $\sqrt 7/2$.
Using this information and~(\ref{ell:eqn}) it is now easy to prove
that the slopes on any component of $\partial N$ having $H$-length
less that $2\pi$ are those in the following set:
\begin{eqnarray*}
\tilde S_1 & =& \Big\{
\infty, -6, -5, -4, -3,
-\frac 52, -2, - \frac 53, -\frac 32, -\frac 43,
-\frac 54, -1, \\ & & \hspace{.5cm}-\frac 34, -\frac 23, -\frac 12,
-\frac 13, -\frac 14, 0, \frac 14, \frac 13,
\frac 12, \frac 23, 1, \frac 32, 2,
3, 4, 5\Big\}.
\end{eqnarray*}

Now we remark that $\tilde S_1\cap\calG=\{\infty,-2,-1,0,1\}$ (recall that the homology
bases are changed: a slope previously corresponding to $p/q$ is now given
by $-2-p/q$, and conversely).
And we set $S_1=\tilde S_1\setminus \calG$, noting that $S_1$ has $23$ elements.
Using SnapPea it is now easy to check that $N(\alpha)$ is actually hyperbolic for all
$\alpha\in S_1$. Moreover, since there is a symmetry of $N$ fixing any boundary component
and exchanging the other two, the two cusps of $N(\alpha)$ have the same shape
given by a certain $z_\alpha\in\matC$. Let us consider as above the maximal cusp
section $H_\alpha$ in $N(\alpha)$ consisting of two tori of the same area, and let us
denote by $A_\alpha$ this area. Again using SnapPea one can determine $z_\alpha$
and $A_\alpha$ for all $\alpha\in S_1$, and correspondingly compute,
via~(\ref{ell:eqn}), the set $\tilde S_2(\alpha)$ of slopes on (any component of)
$\partial N(\alpha)$ having $H_\alpha$-length less than $2\pi$. We now
define
$$S_2(\alpha)=\big\{\beta\in\tilde S_2(\alpha):\ (\alpha,\beta)\not\in\calG\big\},
\qquad S_2=\bigcup\nolimits_{\alpha\in S_1} S_2(\alpha).$$
This set $S_2$ consists of $55$
elements.  Again using SnapPea one can now check that $N(\alpha,\beta)$
is actually hyperbolic for all $(\alpha, \beta)\in (S_1\times
S_2)\setminus\calG$,
and compute the shape $z_{\alpha,\beta}=x_{\alpha,\beta}+i y_{\alpha,\beta}$
and the area $A_{\alpha,\beta}$ of the maximal cusp $H_{\alpha,\beta}$ of
$N(\alpha,\beta)$. With the list of actual values at hand one can check that
$$-1\leqslant x_{\alpha,\beta} \leqslant 4,\quad
0.5 \leqslant y_{\alpha,\beta} \leqslant 4,\quad
A_{\alpha,\beta}\geqslant 1\qquad\forall (\alpha,\beta)\in (S_1\times S_2)\setminus\calG$$
(the last inequality is actually known from the theory, but one can also check it
experimentally).
Now let $p/q$ be a slope on $\partial N(\alpha,\beta)$. Up to switching
signs we can assume $q\geqslant 0$, whence, by~(\ref{ell:eqn}) and the previous
inequalities,
$$\ell(p/q)^2  =  \frac {A_{\alpha,\beta}}{y_{\alpha,\beta}}
\big((p+q\cdot x_{\alpha,\beta})^2+(q\cdot y_{\alpha,\beta})^2\big)
\geqslant \frac 14\big(f(p,q)^2+q^2\cdot (0.5)^2\big) =: g(p,q)$$
where $f(p,q)$ is $0$ if $0\in(p-q,p+4q)$, and $\min\{|p-q|,|p+4q|\}$ otherwise.
We then define $S_3=\{p/q:\ g(p,q)\leqslant 4\pi^2\}$, and
note that $S_3$ has $200$ elements.

Having constructed $S_1,S_2,S_3$, let us check that they work.
If $N(\alpha,\beta,\gamma)$ is not negatively curved then, up to
permutation, we have $\alpha\in\tilde S_1$ by the $2\pi$-theorem.
So either $\alpha\in S_1$ or $\alpha\in\calG$.
And $(\alpha,\beta,\gamma)\in\calG$ in the latter case.
Assuming $\alpha\in S_1$, we know $N(\alpha)$ is hyperbolic, so
$\beta\in\tilde S_2(\alpha)$ up to permutation by the $2\pi$-theorem.
So either $\beta\in S_2$, or $(\alpha,\beta)\in\calG$, or both.
And $(\alpha,\beta,\gamma)\in\calG$ if $(\alpha,\beta)\in\calG$.
Assuming $(\alpha,\beta)\not\in\calG$ we know that $N(\alpha,\beta)$
is hyperbolic, so $\gamma\in S_3$ by the $2\pi$-theorem.

\paragraph{Step~4: Finite enumeration}
To conclude we are left to investigate
all the $23\cdot 55\cdot 200 =253000$
triples of slopes in $S_1\times S_2\times S_3$, and check that those not
lying in $\calG$ define hyperbolic fillings. We have done this using the
scriptable Python version of SnapPea.
Many of the triples in
$(S_1\times S_2\times S_3)\setminus\calG$ turned out to give solutions of the
hyperbolic Dehn filling equations with some inverted tetrahedra.
For these triples we have tried with several alternative
triangulations of $N$ and of some partial fillings of $N$,
getting a genuine solution after some attempts.
The only triple resisting to this method was $\{1,-5,-1/2\}$,
which is known~\cite{Ho-We} to give the closed hyperbolic manifold with
volume 1.0149416\ldots (and
complexity 9).
\finedimo

\appendix

\section{Hyperbolic manifolds with of 1 or 2 cusps}\label{other:fillings:section}
Theorems~\ref{one:cusp:teo},~\ref{two:cusps:teo}, and
\ref{three:cusps:teo} can be used to list all the
non-hyperbolic Dehn fillings of infinitely many hyperbolic manifolds with 1 or 2
cusps. We study in this section the most notable of these manifolds.
Recall that $N$ always denotes the chain-link complement.

\subsection{Number of exceptional slopes} \label{number:subsection}
If $M$ is a 1-cusped hyperbolic manifold, following~\cite{Kirby},
we denote by $E(M)$ the set of exceptional slopes on $\partial M$
(\emph{i.e}~the slopes giving non-hyperbolic filling),  and we set
$e(M)=|E(M)|$, which is finite by~\cite{bibbia}. If $M$ has more
than one cusp, we pick one $T$ and define $E_T(M)$ as the set of
exceptional slopes on $T$. Then we define $e(M)$ as the maximum of
$|E_T(M)|$ over all cusps $T$. Theorem~\ref{one:cusp:teo} readily
implies the following:

\begin{cor}
$E_T(N)=\{\infty,-3,-2,-1,0\}$ for all cusps $T$,
whence $e(N)=5$.
\end{cor}

We do not know if there is any other hyperbolic $M$
with at least 3 cusps and $e(M)\geqslant 5$.

\paragraph{Exceptional slopes on 2-cusped manifolds}
Since no filling of a graph manifold is hyperbolic, Theorem~\ref{one:cusp:teo}
gives $e(M)\geqslant 5$ for all hyperbolic partial fillings $M$ of $N$. We carefully
describe now the $M$'s such that $e(M)>5$, starting from the 2-cusped ones.
Theorems~\ref{one:cusp:teo} and \ref{two:cusps:teo} and
Proposition~\ref{repetitions:prop} imply the following:

\begin{cor}\label{e:Nalpha:cor}
Assume $\alpha\not\in\{\infty,-3,-2,-1,0\}$, so
that $N(\alpha)$ is hyperbolic.  Then $e(N(\alpha))=5$
except for $\alpha=-4,-5/2,-3/2,-1/2,1$, in which case
$e(N(\alpha))=6$. There are $4$ distinct $N(\alpha)$'s
corresponding to these $5$ slopes, because $N(-4) = N(-3/2)$.
\end{cor}

Many informations on the 4 manifolds described in the previous
proposition and their non-hyperbolic fillings are contained in
Table~\ref{four:manifolds:table}. All 4 manifolds are actually
well-known ones, being the Whitehead link exterior, the Whitehead
sister (or $(-2,3,8)$ pretzel) link exterior, the Berge manifold
and the exterior of the 2-bridge link with parameter $3/10$
(Remarks~\ref{Berge:rem} and~\ref{other:2:cusped:rem} below show
how to recognize each of them by looking at its fillings). We
provide for each manifold its label $Mi_j^k$ in the
Callahan-Hildebrand-Weeks census~\cite{Ca-Hi-We} (recall that
$Mi_j^k$ is the $j$-th manifold among those having $k$ cusps and
an ideal triangulation with $i$ tetrahedra, ordered by increasing
volume, and $k$ is omitted if it is $1$). We also show the
relevant slopes, addressing the reader to Fig.~\ref{N:fig} for an
explanation of their pictorial description. Finally, we completely describe
the non-hyperbolic fillings.

An important \emph{caveat} concerning Table~\ref{four:manifolds:table}
is that, to express slopes on $\partial N(\alpha)$,
we are still using  the homology bases induced from those of $N$.
In particular, when $N(\alpha)$ is a link complement, \emph{we are
not using the natural (meridian, longitude) bases}.  In some
cases switching from our basis to the natural one is easy (as an example,
for the Whitehead link complement, using the natural basis each coefficient $\beta$ in
Table~\ref{four:manifolds:table} should be replaced by $\beta-1$),
but in other cases it is harder. For $M4_2^2$ the list of exceptional
slopes refers to its realization as $N(-3/2)$, not as $N(-4)$.

\begin{rem} \label{Berge:rem}
\emph{
The Berge manifold $N(-5/2)$
is the unique manifold having $3$ fillings giving the solid
torus~\cite{Be}. As shown in Table~\ref{four:manifolds:table},
it also has $3$ annular and toroidal fillings.
In~\cite{FriMaPe3} we have generalized this example, building
an infinite class of hyperbolic manifolds $M$ with non-empty geodesic
boundary and $e(M)=6$,
with 3 fillings giving a handlebody and 3 other ones giving an
annular manifold.}
\end{rem}

\begin{rem} \label{other:2:cusped:rem}
\emph{
The manifold $N(-3/2)$ is determined as the only one
with two annular fillings at distance $5$,
while $N(1)$ and $N(-1/2)$ are the only two manifolds
with two annular fillings at distance $4$, see~\cite{gordon, GoWu-TA}.
They can be further distinguished because $N(1)$ has one
filling containing a non-separating torus,
while $N(-3/2)$ has not.
}
\end{rem}

\begin{rem}
\emph{
The $4$ manifolds of Table~\ref{four:manifolds:table}
are the 2-cusped ones with least
known volume, and precisely the $4$ ones having
smallest complexity $4$ (see~\cite{Ca-Hi-We}). The first 2
have canonical decomposition given by one regular ideal
octahedron, the other 2 by $4$ regular ideal tetrahedra.
Repetitions in Table~\ref{four:manifolds:table}
are the effect of the relations listed in Proposition~\ref{repetitions:prop}.}
\end{rem}

\begin{table}
\begin{center}
\begin{tabular}{|c|c|}
\hline

\begin{minipage}{.45\textwidth}
    \phantom{\small{|}}
    \centering
    \mettifig{Whitehead, width = .9\textwidth}
    The Whitehead link exterior $M4_1^2 = N(1)$ \phantom{\Big|}
\end{minipage}
&
\begin{minipage}{.5\textwidth}
  $\begin{array}{rl}
        \infty: & D\times S^1 \phantom{\Big|} \\
        -3: & \seifdue D2121 \bigu 0110\!\! \seifuno A21 \\
        -2: & \seifdue D3131 \\
        -1: & \seifdue D2141 \\
        0: & \seifdue D2131 \\
        1:  & \big(P\times S^1\big) \bigb 0110
    \end{array}$
\end{minipage} \\
\hline

\begin{minipage}{.45\textwidth}
    \phantom{\small{|}}
    \centering
    \mettifig{pretzelm238, width = .9\textwidth}
    The Whitehead sister link exterior $M4_2^2 = N(-3/2) = N(-4)$ \phantom{\Big|}
\end{minipage}
&
\begin{minipage}{.5\textwidth}
  $\begin{array}{rl}
        \infty, -2: & D\times S^1 \phantom{\Big|} \\
        -3, -1: & \seifdue D2131 \\
        -5/2, 0: & \seifdue D2131 \bigu 1110\\
        & \seifuno A21
    \end{array}$
\end{minipage} \\
\hline

\begin{minipage}{.45\textwidth}
    \phantom{\small{|}}
    \centering
  \mettifig{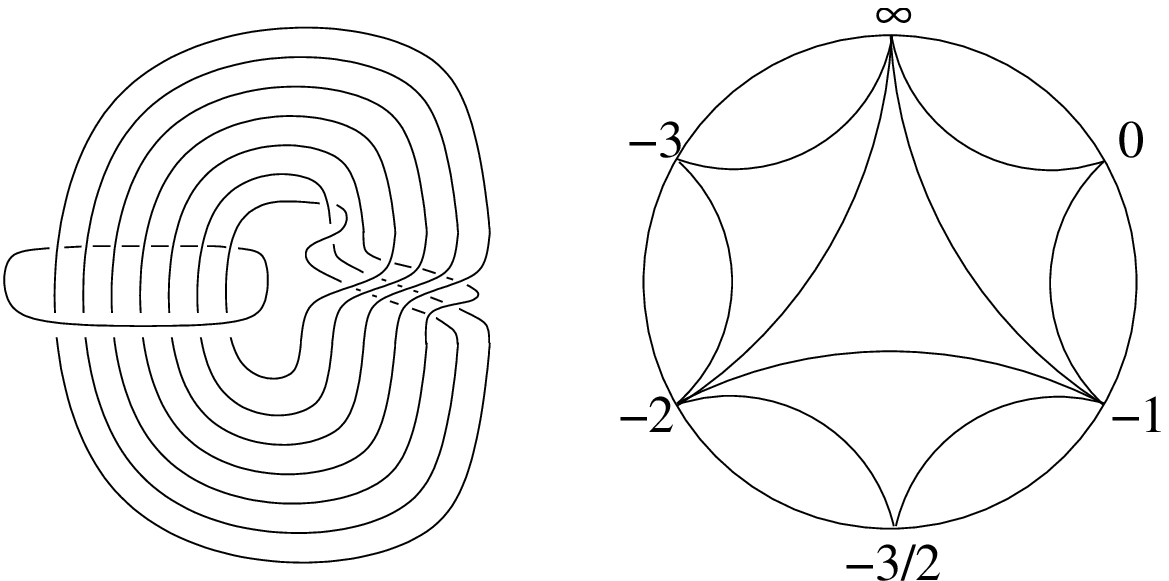, width=.9\textwidth}
    The Berge manifold $M4_3^2 = N(-5/2)$
    \phantom{\Big|}
\end{minipage}
&
\begin{minipage}{.5\textwidth}
  $\begin{array}{rl}
    \infty,-2,-1: & D\times S^1 \phantom{\Big|} \\
        -3,-3/2,0: & \seifdue D2131 \bigu 1110 \\
        & \seifuno A21 \\
    \end{array}$
\end{minipage} \\
\hline

\begin{minipage}{.45\textwidth}
    \phantom{\small{|}}
    \centering
    \mettifig{bridge310, width = .9\textwidth}
    The $3/10$-bridge link exterior $M4_4^2 = N(-1/2)$
    \phantom{\Big|}
\end{minipage}
&
\begin{minipage}{.5\textwidth}
  $\begin{array}{rl}
    \infty: & D\times S^1 \phantom{\Big|} \\
        -4,0: & \seifdue D2131 \bigu 1110\!\! \seifuno A21 \\
        -3,-1: & \seifdue D2152 \\
        -2: & \seifdue D3132
    \end{array}$
\end{minipage} \\
\hline

\end{tabular}
\nota{Non-hyperbolic fillings of the 2-cusped $N(\alpha)$'s with $e(N(\alpha))>5$}
\label{four:manifolds:table}
\end{center}
\end{table}

\paragraph{Exceptional slopes on 1-cusped manifolds}
We now concentrate on the hyperbolic $N(\alpha,\beta)$'s and describe
those having more than 5 exceptional slopes.
As above, we first provide a summarizing statement and then give
many extra informations.

\begin{cor}\label{e:Nalphabeta:cor}
Assume $\alpha,\beta$
do not appear (up to permutation) in
Table~\ref{two:cusps:table}, so $M=N(\alpha,\beta)$ is hyperbolic.
Then:
\begin{description}
\item[(A)] $e(M)\geqslant 7$ for the $11$ manifolds described in Tables~\ref{one:cusp:table1}
to~\ref{one:cusp:table3} (where $e(M)$ attains once the value $10$, twice
the value $8$, and $8$ times the value $7$);
\item[(B)] $e(M)=6$ if $\alpha\in\{-4,-5/2,-3/2,-1/2,1\}$  and $\beta$ is arbitrary,
provided $(\alpha,\beta)$ is not in Tables~\ref{one:cusp:table1}
to~\ref{one:cusp:table3}; the corresponding fillings are as
in Tables~\ref{two:cusps:table1} to~\ref{two:cusps:table4};
in all cases $E(N(\alpha,\beta))=E(N(\alpha))$;
\item[(C)] $e(M)=6$ for the $7$ manifolds described in Table~\ref{one:cusp:table4};
\item[(D)] $e(M)=5$ in all other cases.
\end{description}
No $M$ is of both types (B) and (C); if it is of type (B) then it can be realized
as $N(\alpha,\beta)$ for a unique $\alpha\in\{-5/2,-3/2,-1/2,1\}$.
\end{cor}
\begin{proof}
The only non-trivial fact to be proved is the last sentence. Suppose
$e(M)=6$ and $M=N(\alpha,\beta)$ is of type (B) with
$\alpha\in\{-5/2,-3/2,-1/2,1\}$.
If $M$ is also of type (C), its non-hyperbolic fillings are listed in
Table~\ref{one:cusp:table4}.
It is not difficult to see from Tables~\ref{two:cusps:table1}
to~\ref{two:cusps:table4} that no $M$ of type (B) can have these
fillings, so $M$ is not of type (C). Assume $M=N(\alpha',\beta')$ for
some $\alpha'\in\{-5/2,-3/2,-1/2,1\}$ distinct from $\alpha$.
Since $E(N(\alpha,\beta))$ and $E(N(\alpha',\beta'))$ are the same, we must
have $\{\alpha,\alpha'\}=\{1,-1/2\}$. But $N(1,\beta,1)$ contains a
non-separating torus, whereas $N(-1/2,\beta',\gamma)$ does not,
whatever $\gamma$.
\end{proof}

We warn the reader that, as above, the filling coefficients of each $N(\alpha,\beta)$ in
Tables~\ref{one:cusp:table1} to~\ref{one:cusp:table4} are expressed
\emph{using the homology bases induced from $N$}.
When the same manifold admits several presentations as $N(\alpha,\beta)$,
the coefficients refer to the first presentation listed.
In Tables~\ref{one:cusp:table1} to~\ref{one:cusp:table3},
we also provide for each manifold a surgery presentation along on a link
in $S^3$, taken from~\cite{Kirby}.
Surgery coefficients on link diagrams always
refer to the meridian-longitude homology bases.

Our next remark concerning Corollary~\ref{e:Nalphabeta:cor}
is that the 11 manifolds of point (A) are precisely
the 11 known $M$'s with $e(M)\geqslant 7$, and they are conjectured~\cite{Kirby} to be
the only ones. Moreover:
\begin{itemize}
\item precisely 2 of the 11 manifolds are knot complements,
namely $N2_1=N(1,2)$, the complement of the
figure-8 knot, and $M3_7$ (obtained in 6 different ways as a
filling of $N$), the complement of the $(-2,3,7)$ pretzel knot;
\item all the 11 manifolds except $M3_7$ and $M4_{10}=N(2,2)$
are fillings of the Whitehead link; $M3_7$ is a filling of
any of the $4$ manifolds in Table~\ref{four:manifolds:table}
except the Whitehead link; $M4_{10}$ is not a filling of any
of these manifolds.
\end{itemize}

Turning to the $M$'s with $e(M)=6$
in Corollary~\ref{e:Nalphabeta:cor}, we note that point (B) of course gives
infinitely many examples, thoroughly
described in Tables~\ref{two:cusps:table1} to~\ref{two:cusps:table4}
according to which of the $4$ manifolds in Table~\ref{four:manifolds:table}
they are a filling of. The captions of Tables~\ref{two:cusps:table2}
to~\ref{two:cusps:table4} should be interpreted as shortened versions
of the caption of Table~\ref{two:cusps:table1}.

\begin{table}
\begin{center}
\begin{tabular}{|c|c|}
\hline
\begin{minipage}{.35\textwidth}
    \phantom{\small{|}}
  \centering
\altramettifig{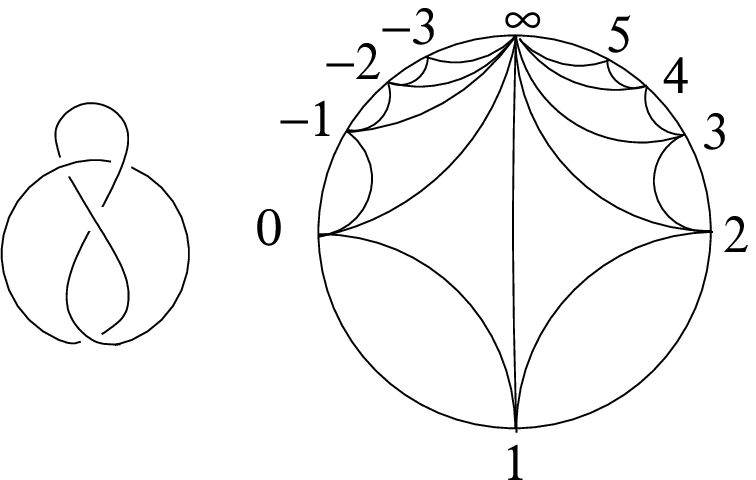}
    The figure-8 knot $M2_1 = N(1,2)$ \phantom{\Big|}
\end{minipage}
&
\begin{minipage}{.60\textwidth}
  $\begin{array}{rl}
        \infty: & S^3 \phantom{\Big|} \\
    -3,5: & \seifdue{D}2121\bigcup\nolimits_{\tiny{\matr0110}}\!\seifdue{D}2131 \\
    -2,4: & \seiftreeul{S^2}313141{-1} \\
    -1,3: & \seiftreeul{S^2}214151{-1} \\
    0,2: & \seiftreeul{S^2}213171{-1} \\
    1:  & T_{\tiny{\matr 31{-1}0}} \\
  \end{array}$ \\
\end{minipage} \\
\hline

\begin{minipage}{.35\textwidth}
    \phantom{\small{|}}
    \centering
\altramettifig{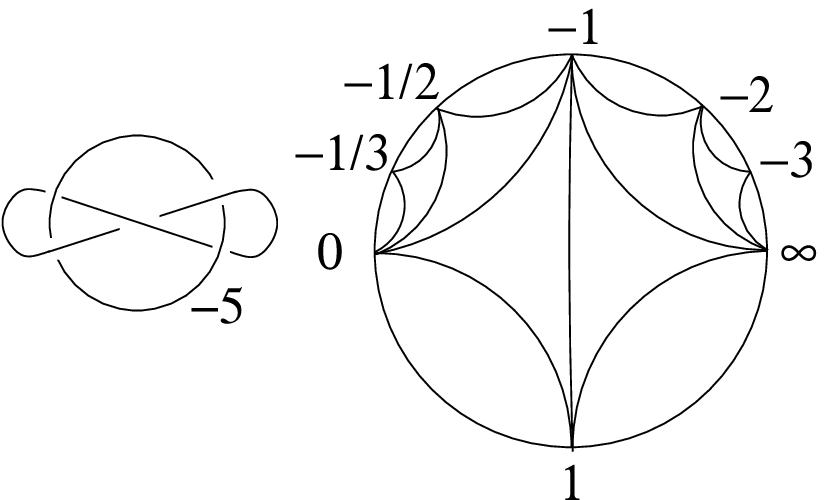}
    The figure-8 knot sister \\
    $M2_2 = N(1,-4) = $ \\
    $N(-4,-3/2) = N(-3/2,-3/2)$ \phantom{\Big|}
\end{minipage}
&
\begin{minipage}{.60\textwidth}
  $\begin{array}{r@{\ }l}
    -1: & \lens{10}3 \\
    -1/2, -2: & \seiftreeul{S^2}213232{-1} \\
    -1/3, -3: & \seifdue{D}2121\!\bigcup\nolimits_{\tiny{\matr01{-1}{-1}}}\!\seifdue{D}2131 \\
    0, \infty: & \lens 51 \\
    1:  & T_{\tiny{\matr {-3}1{-1}0}}
  \end{array}$
\end{minipage} \\
\hline

\begin{minipage}{.35\textwidth}
  \centering
\altramettifig{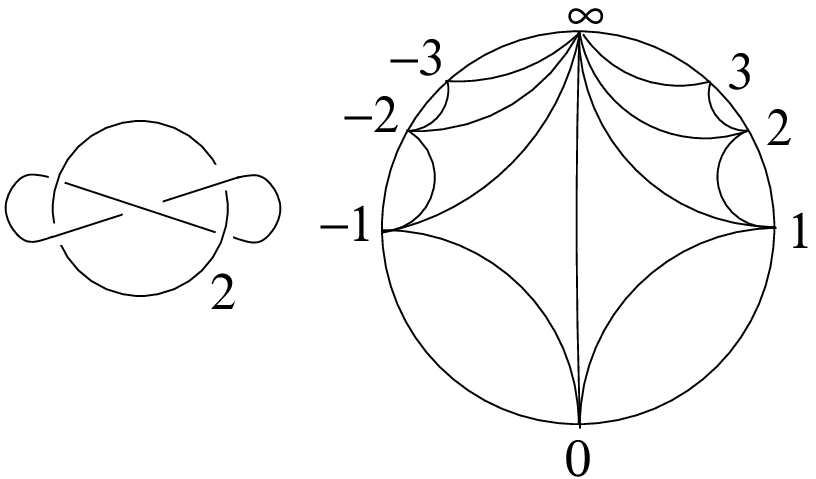}
    $M3_3 = N(1,3)$ \phantom{\Big|}
\end{minipage}
&
\begin{minipage}{.60\textwidth}
  $\begin{array}{r@{\ }l}
    \infty: & \mathbb{P}^3 \phantom{\Big|} \\
    -3: & \seifdue{D}2121\!\bigcup\nolimits_{\tiny{\matr0110}}\!\seifdue{D}2141 \\
    -2: & \seiftreeul{S^2}313151{-1} \\
    -1: & \seiftreeul{S^2}214161{-1} \\
    0:  & \seiftreeul{S^2}213181{-1} \\
    1:  & T_{\tiny{\matr 41{-1}0}} \\
    2:  & \seiftreeul{S^2}214151{-1} \\
    3: & \seifdue{D}2121\!\bigcup\nolimits_{\tiny{\matr120{-1}}}\!\seifdue{D}2131
    \phantom{|}_{\phantom{\Big|}}\\
  \end{array}$
\end{minipage} \\
\hline

\end{tabular}
\nota{Non-hyperbolic fillings of 1-cusped $M=N(\alpha,\beta)$'s with $e(M)\geqslant 8$} \label{one:cusp:table1}
\end{center}
\end{table}

\begin{table}
\begin{center}
\begin{tabular}{|c|c|}
\hline
\begin{minipage}{.35\textwidth}
    \phantom{\small{|}}
  \centering
\altramettifig{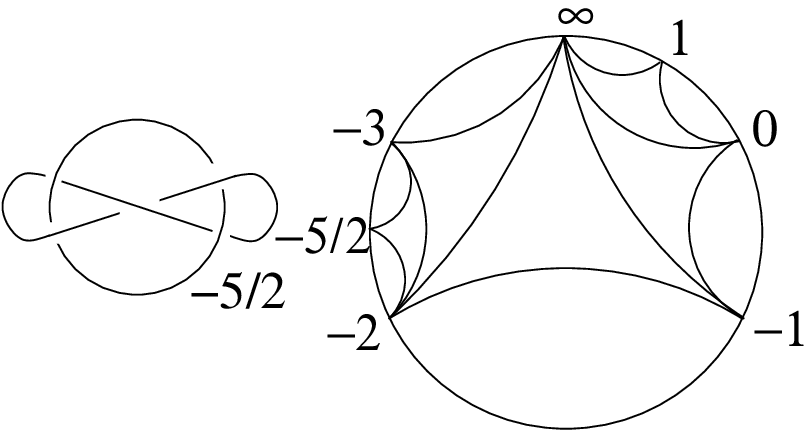}
    $M3_1 = N(1,-3/2) =$ \\
    $N(-4,-4)  = N(-4,-2/3)$ \\
    $ = N(-3/2,-7/3)$ \phantom{\Big|}
\end{minipage}
&
\begin{minipage}{.60\textwidth}
  $\begin{array}{r@{\ }l}
    \infty: & \lens 52 \phantom{\Big|}\\
    -3:   & \seiftre{S^2}212132 \\
    -2:   & \lens {15}4 \\
    -1:   & \seiftreeul{S^2}213241{-1} \\
    0:    & \seiftreeul{S^2}213172{-1} \\
    1:    & \seifuno A21 /_{\tiny{\matr 1110}} \phantom{\big|}\\
    -5/2: & \seifdue{D}2131\!\bigcup\nolimits_{\tiny{\matr110{-1}}}\!\seifdue{D}2131
    \phantom{|}_{\phantom{\Big|}}\\
  \end{array}$
\end{minipage} \\
\hline

\begin{minipage}{.35\textwidth}
    \phantom{\small{|}}
  \centering
\altramettifig{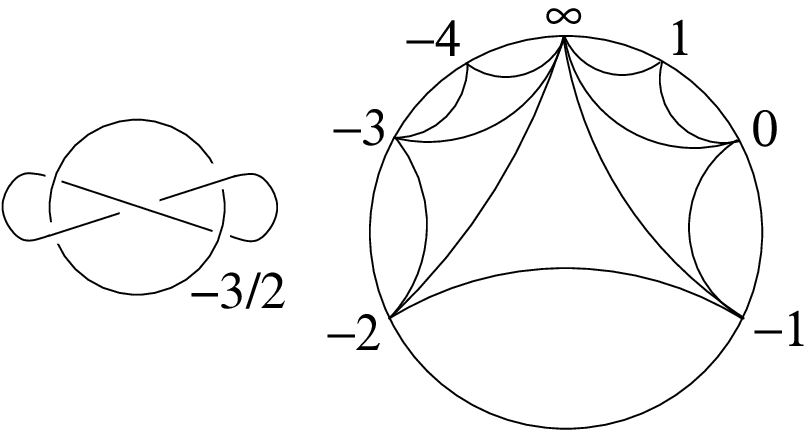}
    $M3_2 = N(1,-1/2) = N(-5,-1/2)$ \phantom{\Big|}
\end{minipage}
&
\begin{minipage}{.60\textwidth}
  $\begin{array}{r@{\ }l}
    \infty: & \lens 31 \phantom{\Big|}\\
    -4:   & \seiftreeul{S^2}213232{-1} \\
    -3:   & \seiftreeul{S^2}212153{-1} \\
    -2:   & \seiftreeul{S^2}313132{-1} \\
    -1:    & \seiftreeul{S^2}214152{-1} \\
    0:    & \seiftreeul{S^2}213192{-1} \\
    1:    & \seifuno A21 /_{\tiny{\matr 0110}}
  \end{array}$
\end{minipage} \\
\hline

\begin{minipage}{.35\textwidth}
    \phantom{\small{|}}
  \centering
\altramettifig{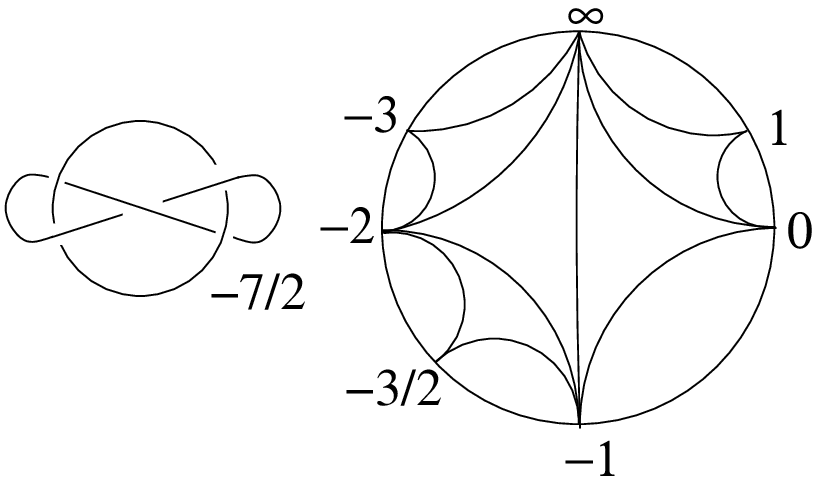}
    $M3_8 = N(1,-5/2) = N(-5/2, -4/3) = N(-5/2, -5/2)$ \phantom{\Big|}
\end{minipage}
&
\begin{minipage}{.60\textwidth}
  $\begin{array}{r@{\ }l}
    \infty: & \lens 72 \phantom{\Big|} \\
    -3:   & \seifdue{D}2121\!\bigcup\nolimits_{\tiny{\matr011{-1}}}\!\seifdue{D}2131 \\
    -2:   & \lens {21}8 \\
    -3/2: & \seifdue{D}2131\!\bigcup\nolimits_{\tiny{\matr110{-1}}}\!\seifdue{D}2131 \\
    -1:   & \lens {14}3 \\
    0:    & \seiftreeul{S^2}213152{-1} \\
    1:    & \seifuno A21 /_{\tiny{\matr 2110}} \\
  \end{array}$
\end{minipage} \\
\hline

\begin{minipage}{.35\textwidth}
    \phantom{\small{|}}
  \centering
\altramettifig{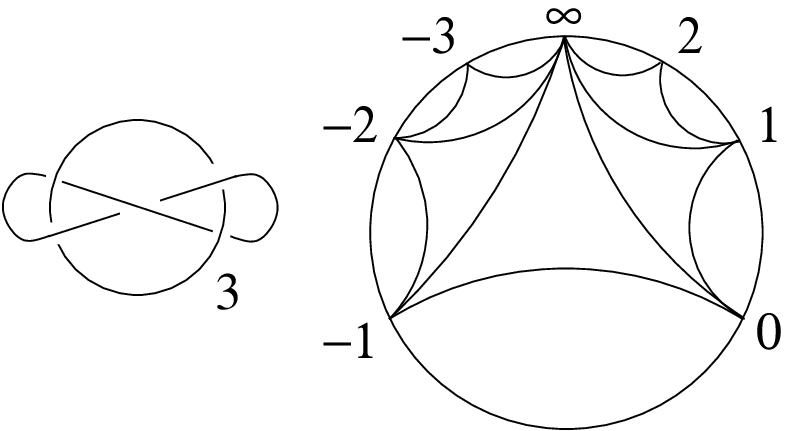}
    $M4_2 = N(1,4)$ \phantom{\Big|}
\end{minipage}
&
\begin{minipage}{.60\textwidth}
  $\begin{array}{r@{\ }l}
    \infty: & \lens 31 \phantom{\Big|}\\
    -3:   & \seifdue{D}2121\!\bigcup\nolimits_{\tiny{\matr0110}}\!\seifdue{D}2151 \\
    -2:   & \seiftreeul{S^2}313161{-1} \\
    -1:   & \seiftreeul{S^2}214171{-1} \\
    0:    & \seiftreeul{S^2}213191{-1} \\
    1:    & T_{\tiny{\matr 51{-1}0}} \\
    2:    & \seiftreeul{S^2}313141{-1}
    \phantom{\Big|}\\
  \end{array}$
\end{minipage} \\
\hline

\end{tabular}
\nota{Non-hyperbolic fillings of 1-cusped $M=N(\alpha,\beta)$'s with $e(M)=7$} \label{one:cusp:table2}
\end{center}
\end{table}

\begin{table}
\begin{center}
\begin{tabular}{|c|c|}
\hline
\begin{minipage}{.35\textwidth}
    \phantom{\small{|}}
  \centering
\altramettifig{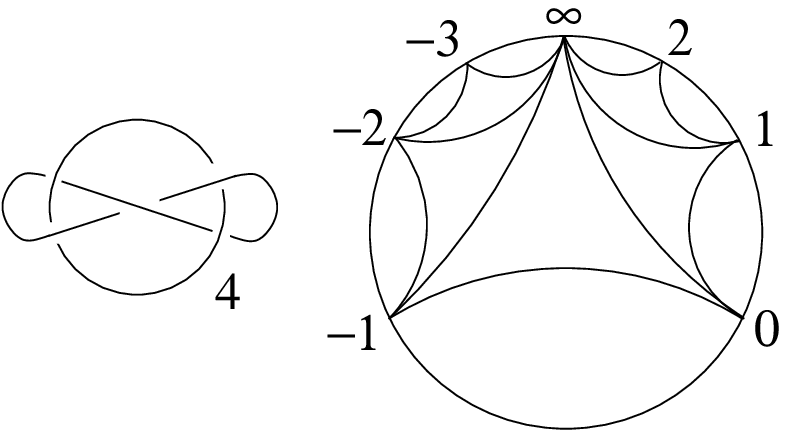}
    $M5_1 = N(1,5)$ \phantom{\Big|}
\end{minipage}
&
\begin{minipage}{.60\textwidth}
  $\begin{array}{r@{\ }l}
    \infty: & \lens 41 \phantom{\Big|}\\
    -3:   & \seifdue{D}2121\!\bigcup\nolimits_{\tiny{\matr0110}}\!\seifdue{D}2161 \\
    -2:   & \seiftreeul{S^2}313171{-1} \\
    -1:   & \seiftreeul{S^2}214181{-1} \\
    0:    & \seiftreeul{S^2}2131{10}1{-1} \\
    1:    & T_{\tiny{\matr 61{-1}0}} \\
    2:    &
    \seifdue{D}2121\!\bigcup\nolimits_{\tiny{\matr0110}}\!\seifdue{D}2131
    \phantom{|}_{\phantom{\Big|}}\\
  \end{array}$
\end{minipage} \\
\hline

\begin{minipage}{.35\textwidth}
    \phantom{\small{|}}
  \centering
\altramettifig{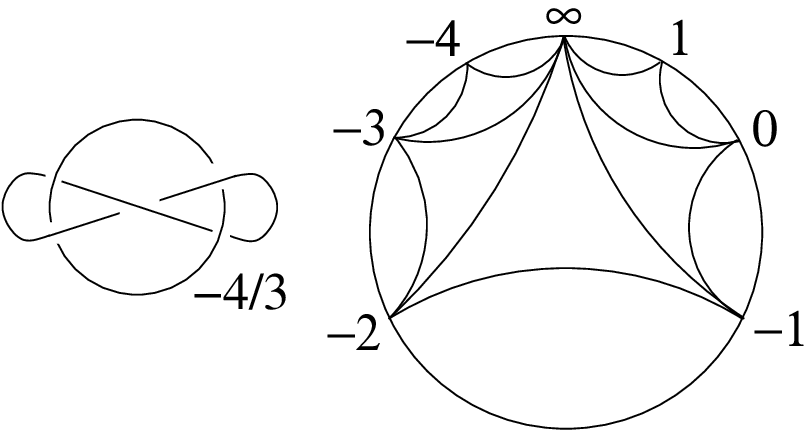}
    $M4_{12} = N(1,-1/3)$ \phantom{\Big|}
\end{minipage}
&
\begin{minipage}{.60\textwidth}
  $\begin{array}{r@{\ }l}
    \infty: & \lens 41 \phantom{\Big|}\\
    -4:   & \seifdue{D}2121\!\bigcup\nolimits_{\tiny{\matr01{-1}{-1}}}\!\seifdue{D}2131 \\
    -3:   & \seifdue{D}2121\!\bigcup\nolimits_{\tiny{\matr0111}}\!\seifdue{D}2121 \\
    -2:   & \seiftreeul{S^2}313153{-1} \\
    -1:   & \seiftreeul{S^2}214183{-1} \\
    0:   & \seiftreeul{S^2}2131{14}3{-1} \\
    1:    & \seifuno A32 /_{\tiny{\matr 0110}}
    \phantom{|}_{\phantom{\Big|}}\\
  \end{array}$
\end{minipage} \\
\hline

\begin{minipage}{.35\textwidth}
    \phantom{\small{|}}
  \centering
\altramettifig{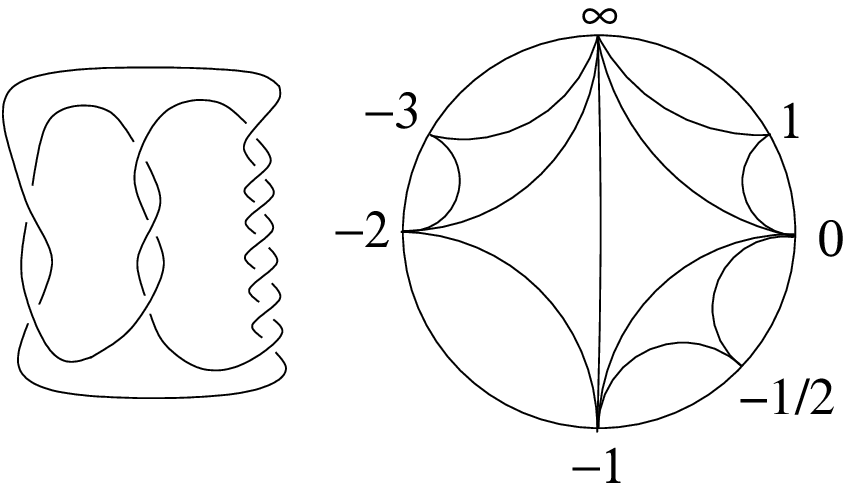}
    $(-2,3,7)$ pretzel knot
    ${\scriptstyle M3_7 = N(-4,-1/3) = N(-3/2,-1/2)}$ \\
    ${\scriptstyle = N(-3/2,-8/3) = N(-4,-5/2)}$ \\
    ${\scriptstyle = N(-1/2,-5/2) = N(-5/2,-5/3)}$
     \phantom{\Big|}
\end{minipage}
&
\begin{minipage}{.60\textwidth}
  $\begin{array}{r@{\ }l}
    \infty: & S^3 \phantom{\Big|}\\
    -3:   & \seifdue{D}2121\!\bigcup\nolimits_{\tiny{\matr{-1}10{-1}}}\!\seifdue{D}2131 \\
    -2:   & \seiftreeul{S^2}213252{-1} \\
    -1:   & \lens{18}5 \\
    -1/2: & \seifdue{D}2131\!\bigcup\nolimits_{\tiny{\matr11{-2}{-1}}}\!\seifdue{D}2131 \\
    0:    & \lens{19}7 \\
    1:    & \seifdue{D}2121\!\bigcup\nolimits_{\tiny{\matr01{-1}{-1}}}\!\seifdue{D}2131 \\
  \end{array}$
\end{minipage} \\
\hline

\begin{minipage}{.35\textwidth}
    \phantom{\small{|}}
  \centering
\altramettifig{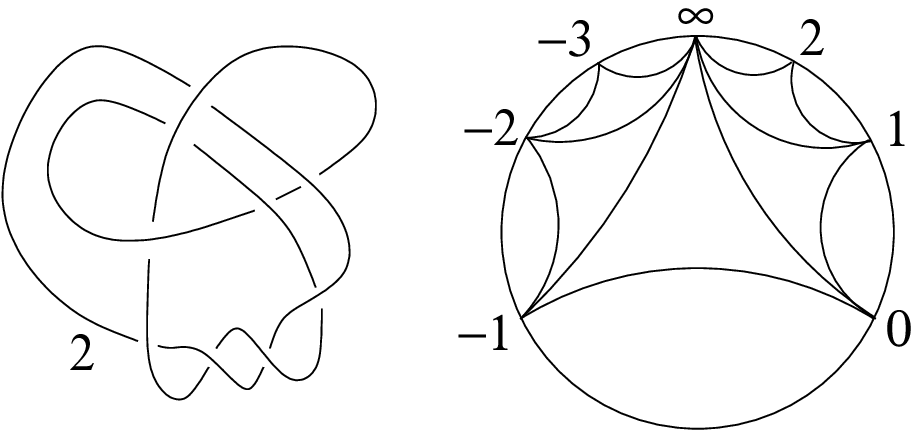}
    $M4_{10} = N(2,2)$ \phantom{\Big|}
\end{minipage}
&
\begin{minipage}{.60\textwidth}
  $\begin{array}{r@{\ }l}
    \infty: & \lens 31 \phantom{\Big|}\\
    -3:   & \seifdue{D}2131\!\bigcup\nolimits_{\tiny{\matr0110}}\!\seifdue{D}2131 \\
    -2:   & \seiftreeul{S^2}314141{-1} \\
    -1:   & \seiftreeul{S^2}215151{-1} \\
    0:    & \seiftreeul{S^2}213181{-1} \\
    1:    & \seiftreeul{S^2}213171{-1} \\
    2:    & \seifdue{D}2121\!\bigcup\nolimits_{\tiny{\matr23{-1}{-2}}}\!\seifdue{D}2131
    \phantom{|}_{\phantom{\Big|}}\\
  \end{array}$
\end{minipage} \\
\hline

\end{tabular}
\nota{Non-hyperbolic fillings of 1-cusped $M=N(\alpha,\beta)$'s with $e(M)=7$ (continued)} \label{one:cusp:table3}
\end{center}
\end{table}

\begin{table}
\begin{center}
\begin{tabular}{c|c|l}
$p/q$ & $r/s$ & \qquad $M$ \\
        \hline\hline
$\infty$ & any & $\lens {r-s}s$ \phantom{\Big|} \\
        \hline
$-3$ &
    $-1+\frac 1n$ \phantom{\Big|} &
        $\seiftre{S^2}2121{2n+1}{-2}$ \\
        \cline{2-3}
     &
    $\neq -1+\frac 1n$ \phantom{\Big|} &
    $\seifdue{D}2121\!\bigcup\nolimits_{\tiny{\matr0110}}\!\seifdue{D}21{r+s}s$ \\
        \hline
$-2$ &
    $-2+\frac 1n$ \phantom{\Big|} &
    $\lens{9n-3}{3n-2}$ \\
        \cline{2-3}
     &
    $\neq -2+\frac 1n$ \phantom{\Big|} &
    $\seiftre{S^2}3{-2}31{r+2s}s$ \\
        \hline
$-1$ &
    $-3+\frac 1n$ \phantom{\Big|} &
    $\lens{8n-2}{4n+1}$ \\
        \cline{2-3}
     &
    $\neq -3+\frac 1n$ \phantom{\Big|} &
    $\seiftre{S^2}2{-1}41{r+3s}s$ \\
        \hline
$0$ &
    $-5$ \phantom{\Big|} &
    $\matRP^3 \# \lens 31$ \\
        \cline{2-3}
        &
    $-5+\frac 1n$ \phantom{\Big|} &
    $\lens{6n-1}{2n-1}$ \\
        \cline{2-3}
     &
    $\neq -5, -5+\frac 1n$ \phantom{\Big|} &
    $\seiftre{S^2}2{-1}31{r+5s}s$ \\
        \hline
$1$  &
    $n$ \phantom{\Big|} &
    $T_{\tiny{\matr {n+1}1{-1}0}}$ \\
        \cline{2-3}
     &
    $\not\in\matZ$ \phantom{\Big|} &
    $\seifuno As{r+s} /_{\tiny{\matr 0110}}$ \\
\end{tabular}
\end{center}
\nota{If $M4_1^2=N(1)$ is the Whitehead link exterior and $M=M4_1^2(r/s)=N(1,r/s)$
then $M$ is hyperbolic with $e(M)=6$ if and only if
$r/s\not\in\{\infty,-4,-3,$ $-5/2,-2,-3/2, -1,-1/2,-1/3,$ $0,1/2,0,1,2,3,4,5\}$.
Under this assumption we have $E(M)=\{\infty,-3,-2,-1,0,1\}$
in all cases, and the corresponding non-hyperbolic fillings are as described} \label{two:cusps:table1}
\end{table}

\begin{table}
\begin{center}
\begin{tabular}{c|c|l}
$p/q$ & $r/s$ & \qquad $M$ \\
        \hline\hline
$\infty$ & any & $\lens {3r+2s}{r+s}$ \phantom{\Big|} \\
        \hline
$-3$ &
    $-1+\frac 1n$ \phantom{\Big|} &
    $\lens{6n+7}{3n+2}$ \\
        \cline{2-3}
     &
    $\neq -1+\frac 1n$ \phantom{\Big|} &
    $\seiftre{S^2}2132{r+s}s$ \\
        \hline
$-5/2$ &
    $-2+\frac 1n$ \phantom{\Big|} &
    $\seiftre{S^2}213{-1}{2n-1}2$ \\
        \cline{2-3}
     &
    $\neq -2+\frac 1n$ \phantom{\Big|} &
    $\seifdue{D}2131\!\bigcup\nolimits_{\tiny{\matr110{-1}}}\!\seifdue{D}21{r+2s}s$ \\
        \hline
$-2$ & any & $\lens{4r+11s}{r+3s}$ \phantom{\Big|} \\
        \hline
$-1$  &
    $-3+\frac1n$ \phantom{\Big|} &
    $\lens{6n+1}{3n+2}$ \\
        \cline{2-3}
     &
    $\neq -3+\frac1n$ \phantom{\Big|} &
    $\seiftre{S^2}2{-1}32{r+3s}s$ \\
        \hline
$0$ &
    $n$ \phantom{\Big|} &
        $\seiftre{S^2}2{-1}31{2n+5}2$ \\
        \cline{2-3}
     &
    $\not\in\matZ$ \phantom{\Big|} &
    $\seifdue{D}2131\!\bigcup\nolimits_{\tiny{\matr11{-1}0}}\!\seifdue{D}21s{r+2s}$ \\
\end{tabular}
\end{center}
\nota{Fillings of $M4_2^2=N(-3/2)$ for
$r/s\not\in\{\infty,-4,-3,-8/3,-5/2,-7/3,-2,-3/2,-1,$ $-1/2,0,1\}$}
\label{two:cusps:table2}
\end{table}

\begin{table}
\begin{center}
\begin{tabular}{c|c|l}
$p/q$ & $r/s$ & \qquad $M$ \\
        \hline\hline
$\infty$ & any & $\lens {5r+2s}{2r+s}$ \phantom{\Big|} \\
        \hline
$-3$ &
    $-1+\frac1n$ \phantom{\Big|} &
    $\seiftre{S^2}2{-1}32{2n+1}2$ \\
        \cline{2-3}
     &
    $\neq -1+\frac1n$ \phantom{\Big|} &
    $\seifdue{D}2131\!\bigcup\nolimits_{\tiny{\matr1110}}\!\seifdue{D}21{r+s}s$ \\
        \hline
$-2$ & any & $\lens{8r+13s}{3r+5s}$ \phantom{\Big|} \\
        \hline
$-3/2$ &
    $-2+\frac1n$ \phantom{\Big|} &
    $\seiftre{S^2}2{-1}32{2n-1}2$ \\
        \cline{2-3}
     &
    $\neq-2+\frac1n$ \phantom{\Big|} &
    $\seifdue{D}2131\!\bigcup\nolimits_{\tiny{\matr110{-1}}}\!\seifdue{D}21{r+2s}s$ \\
        \hline
$-1$ & any & $\lens{3r+11s}{2r+7s}$ \phantom{\Big|} \\
        \hline
$0$  &
    $n$ \phantom{\Big|} &
    $\seiftre{S^2}2{-1}31{2n+3}2$ \\
        \cline{2-3}
     &
    $\not\in\matZ$ \phantom{\Big|} &
    $\seifdue{D}2131\!\bigcup\nolimits_{\tiny{\matr1101}}\!\seifdue{D}21sr$ \\
\end{tabular}
\end{center}
\nota{Fillings of $M4_3^2=N(-5/2)$ for
$r/s\not\in\{\infty,-4,-3,-5/2,-2,-5/3,-3/2,-4/3,-1,$ $-1/2,0,1\}$}
\label{two:cusps:table3}
\end{table}

\begin{table}
\begin{center}
\begin{tabular}{c|c|l}
$p/q$ & $r/s$ & \qquad $M$ \\
        \hline\hline
$\infty$ & any & $\lens {r+2s}s$ \phantom{\Big|} \\
        \hline
$-4$ &
    $-\frac12$ \phantom{\Big|}&
    $\seifdueeul{\matRP^2}2131{-1}$ \\
        \cline{2-3}
     &
    $n$ \phantom{\Big|} &
    $\seiftre{S^2}2{-1}32{2n+1}2$ \\
        \cline{2-3}
     &
    $\not\in\matZ\cup\{-\frac12\}$ \phantom{\Big|} &
    $\seifdue{D}2131\!\bigcup\nolimits_{\tiny{\matr1110}}\!\seifdue{D}21sr$ \\
        \hline
$-3$ &
    $-1+\frac1n$ \phantom{\Big|} &
    $\lens{10n+1}{5n-2}$ \\
        \cline{2-3}
     &
    $\neq -1+\frac 1n$ \phantom{\Big|} &
        $\seiftre{S^2}2{-1}53{r+s}s$ \\
        \hline
$-2$ &
    $-2+\frac1n$ \phantom{\Big|} &
    $\lens{9n}{3n-1}$ \\
        \cline{2-3}
     &
    $\neq-2+\frac1n$ &
    $\seiftre{S^2}3{-2}32{r+2s}s$ \phantom{\Big|} \\
        \hline
$-1$ &
    $-3+\frac1n$ \phantom{\Big|} &
    $\lens{10n-1}{5n+2}$ \\
        \cline{2-3}
     &
    $\neq-3+\frac1n$ &
    $\seiftre{S^2}2{-1}52{r+3s}s$ \phantom{\Big|} \\
        \hline
$0$  &
    $-\frac72$ \phantom{\Big|} &
    $\seifdueeul{\matRP^2}2131{-1}$ \\
        \cline{2-3}
     &
    $n$ \phantom{\Big|} &
    $\seiftre{S^2}2{-1}31{2n+7}2$ \\
        \cline{2-3}
     &
    $\not\in\matZ\cup\{-\frac72\}$ \phantom{\Big|} &
    $\seifdue{D}2131\!\bigcup\nolimits_{\tiny{\matr1123}}\!\seifdue{D}21sr$ \\
\end{tabular}
\end{center}
\nota{Fillings of $M4_4^2=N(-1/2)$ for
$r/s\not\in\{\infty,-5,-4,-3,-5/2,-2,-3/2,-1,0,1\}$}
\label{two:cusps:table4}
\end{table}

\begin{table}
\begin{center}
\begin{tabular}{|c|}
\hline
$m139 =
N(-5,-5)$ \phantom{\Big|} \\
$\begin{array}{r@{\ }lr@{\ }l}
    \infty: & \lens{24}5 &
  -3:   & \seifdue{D}2141\!\bigcup\nolimits_{\tiny{\matr110{-1}}}\!\seifdue{D}2141 \\
  -2:   & \seiftre{S^2}313132 &
  -1:   & \seiftre{S^2}212121 \\
    -1/2: & \seifuno A21 /_{\tiny{\matr 0110}} &
  0:    & \seiftreeul{S^2}213261{-1} \\
\end{array}$ \\
\hline
$M4_{35} = N(-5/3,-5/3)$ \phantom{\Big|} \\
$\begin{array}{r@{\ }lr@{\ }l}
    \infty: & \lens{16}7 &
  -5/2: & \seifdue{D}2121\!\bigcup\nolimits_{\tiny{\matr{-1}10{-1}}}\!\seifdue{D}2131 \\
    -3:   & (K,1) &
  -2:   & \lens{16}5 \\
    -1:   & \seiftreeul{S^2}214343{-1} &
  0:   & \seifdue{D}2131\!\bigcup\nolimits_{\tiny{\matr11{-1}0}}\!\seifdue{D}3131 \\
\end{array}$ \\
\hline
$m208 =
N(-7/3,-7/3)$ \phantom{\Big|} \\
$\begin{array}{r@{\ }lr@{\ }l}
    \infty: & \lens{40}{11} &
    -3:   & \seifdue{D}2141\!\bigcup\nolimits_{\tiny{\matr110{-1}}}\!\seifdue{D}2141 \\
  -2:   & \lens{20}7 &
  -3/2: & \seifuno A21 /_{\tiny{\matr 1110}} \\
    -1:   & \seiftreeul{S^2}2121211 &
  0:   & \seifdue{D}2131\!\bigcup\nolimits_{\tiny{\matr1110}}\!\seifdue{D}3131 \\
\end{array}$ \\
\hline
$m120 =
N(2,3)$ \phantom{\Big|} \\
$\begin{array}{r@{\ }lr@{\ }l}
    \infty: & \lens52 &
  -3:   & \seifdue{D}2131\!\bigcup\nolimits_{\tiny{\matr0110}}\!\seifdue{D}2141 \\
    -2:   & \seiftreeul{S^2}314151{-1} &
    -1:   & \seiftreeul{S^2}215161{-1} \\
    0:    & \seiftreeul{S^2}213191{-1} &
    1:    & \seiftreeul{S^2}214151{-1} \\
\end{array}$ \\
\hline
$s90 =
N(2,4)$ \phantom{\Big|} \\
$\begin{array}{r@{\ }lr@{\ }l}
    \infty: & \lens72 &
  -3:   & \seifdue{D}2131\!\bigcup\nolimits_{\tiny{\matr0110}}\!\seifdue{D}2151 \\
    -2:   & \seiftreeul{S^2}314161{-1} &
    -1:   & \seiftreeul{S^2}215171{-1} \\
    0:    & \seiftreeul{S^2}2131{10}1{-1} &
    1:    & \seiftreeul{S^2}313141{-1} \\
\end{array}$ \\
\hline
$v224 =
N(2,5)$ \phantom{\Big|} \\
$\begin{array}{r@{\ }lr@{\ }l}
    \infty: & \lens92 &
  -3:   & \seifdue{D}2131\!\bigcup\nolimits_{\tiny{\matr0110}}\!\seifdue{D}2161 \\
    -2:   & \seiftreeul{S^2}314171{-1} &
    -1:   & \seiftreeul{S^2}215181{-1} \\
    0:    & \seiftreeul{S^2}2131{11}1{-1} &
    1:    & \seifdue{D}2121\!\bigcup\nolimits_{\tiny{\matr0110}}\!\seifdue{D}2131 \\
\end{array}$ \\
\hline
$s118 =
N(3,3)$ \phantom{\Big|} \\
$\begin{array}{r@{\ }lr@{\ }l}
    \infty: & \lens83 &
  -3:   & \seifdue{D}2141\!\bigcup\nolimits_{\tiny{\matr0110}}\!\seifdue{D}2141 \\
    -2:   & \seiftreeul{S^2}315151{-1} &
    -1:   & \seiftreeul{S^2}216161{-1} \\
    0:    & \seiftreeul{S^2}2131{10}1{-1} &
    1:    & \seifdue{D}2121\!\bigcup\nolimits_{\tiny{\matr120{-1}}}\!\seifdue{D}2131 \\
\end{array}$ \\
\hline
\end{tabular}
\end{center}
\nota{Non-hyperbolic fillings of sporadic $M=N(\alpha,\beta)$'s with
$e(M)=6$. The old SnapPea name of $M$ is shown when $M$ has complexity
more than $4$}
\label{one:cusp:table4}
\end{table}

\begin{rem}
{\em
There are $67$ hyperbolic manifolds with complexity at most $4$,
see~\cite{Ca-Hi-We}. Using SnapPea one sees that $60$ of them
are fillings of $N$, while the other $7$, namely
$M4_i$ for $i\in\{33,34,40,41,45,51,52\}$, do not seem to be.
Since $M4_i$ is amphichiral for
$i\in\{33,34,51,52\}$, Proposition~\ref{amphichiral:prop} shows
it is certainly not a filling of $N$ in these cases,
but we do not know for $i\in\{40,41,45\}$.
}
\end{rem}

\paragraph{Automorphisms}
We can now prove Theorem~\ref{cosmetic:teo} and
Proposition~\ref{amphichiral:prop}. We start with the following:

\begin{teo} \label{automorphisms:teo}
Let $M$ be a hyperbolic $N(\alpha)$ or $N(\alpha,\beta)$
and let $f$ be an automorphism of $M$ which leaves
a component $T$ of $\partial M$ invariant. Then $f$ acts
on $H_1(T)$ as the multiplication by $\pm 1$, except when
$M$ and $f$ are as described in Proposition~\ref{repetitions:prop}.
\end{teo}
\begin{proof}
We begin by noting that $f$ acts on the slopes on $T$, and it leaves $E_T(M)$ invariant.
Since $E_T(M)$ contains more than $2$ slopes, the action
on $E_T(M)$ is trivial if and only if the action on
$H_1(T)$ is the multiplication by $\pm 1$.
We suppose now that the action of $f$ on $E_T(M)$ is not trivial and,
examining the various possibilities for $E_T(M)$ and the number
of cusps of $M$, we deduce either a contradiction or that the
action of $f$ is as described in Proposition~\ref{repetitions:prop}.

To begin, we note that the action of $f$ on
$E_T(M)$ must preserve the filled manifolds and can
have at most 2 fixed points. Using also the fact that
$|E_T(M)|\geqslant 5$, we deduce that, among the
non-hyperbolic fillings of $M$, either there is at least one appearing at least
three times, or there are at least two appearing at least twice each.
If $e(M)>6$, Tables~\ref{one:cusp:table1} to~\ref{one:cusp:table3} show that
this happens only for $N(1,2)$ and $N(1,-4)$, and the action of $f$ is
given by relations (\ref{514:eqn}) and (\ref{515:eqn}) respectively.
If $e(M)=6$ and $M$ has two cusps, then $M$ is one of the $4$
manifolds listed in Table~\ref{four:manifolds:table},
and the action of $f$ is given by one of the
relations (\ref{512:eqn}) to (\ref{513:eqn}).

Assume now $e(M)=6$ and $M$ has one cusp.
By Corollary~\ref{e:Nalphabeta:cor}, $M$ is one of the manifolds listed in
Tables~\ref{two:cusps:table1} to~\ref{one:cusp:table4}.
It cannot be one of those in Table~\ref{one:cusp:table4} because in this
table the same filling never appears more than once. The other cases are discussed as follows:
\begin{itemize}
\item If $M$ appears in Table~\ref{two:cusps:table1}
or in Table~\ref{two:cusps:table4}, we have $E(M)=\{\infty,a-2,a-1,a,a+1,a+2\}$
for $a=-1$ or $a=-2$ respectively. Since $\infty$ is the only slope
with distance $1$ from all other ones in $E(M)$, it is fixed by $f$.
The action is then given by $a+n\mapsto a-n$,
but one sees quite easily that $M(a-2)$ is never equal to $M(a+2)$:
a contradiction;
\item If $M$ appears in Table~\ref{two:cusps:table2}, one sees that $f$
is the unique non-trivial symmetry of the set
$E(M)=\{\infty, -3,-5/2,-2,-1,0\}$, which interchanges $\infty$ and $-2$.
But $\lens{3r+2s}{r+s}$ never equals $\lens{4r+11s}{r+3s}$: a contradiction;
\item If $M$ appears in Table~\ref{two:cusps:table3}, one sees that $f$
permutes
non-trivially each set $\{-2,-1,\infty\}$ and $\{-3,-3/2,0\}$.
The three lens spaces $M(\infty), M(-2),$ and $M(-1)$ are distinct
except for $r/s\in\{-2/5,-11/3,-13/8\}$, but in each of these cases
(which give the cosmetic example found in~\cite{cosmetic} and mentioned in
Subsection~\ref{consequences:subsection}) the
manifolds $M(-3), M(-3/2),$ and $M(0)$ are distinct, giving a contradiction.
\end{itemize}

We are left to discuss the case where $e(M)=5$, so $E_T(M)=
\{\infty,-3,-2,-1,0\}$, and $f$ is the unique non-trivial symmetry of
$E_T(M)$, which fixes $\infty$ and maps
$n$ to $-3-n$. With some effort, one checks from
Tables~\ref{two:cusps:table} and~\ref{closed:fill:table1} that $M(0)$
and $M(3)$ are always distinct.
\end{proof}

\begin{rem}
\emph{If $M=N(\alpha,\beta)$ then the assumption $f(T)=T$ is empty.
If $M=N(\alpha)$ there exists $\tau:M\to M$ switching the components
of $\partial M$, so either $f$ or $\tau\circ f$ satisfies the assumption.}
\end{rem}

\dimo{amphichiral:prop}
Let $M$ be amphichiral and let $\tau :M\to M$ be orientation-reversing.
Since the components of $\partial M$ are symmetric, we can assume
$\tau$ leaves them invariant. If $T\subset\partial M$ then
$\tau$ acts on $H_1(T)$ as a matrix with determinant $-1$, in
particular not as multiplication by $\pm 1$. Then $M$ and $\tau$ appear
in Proposition~\ref{repetitions:prop}. And the only
orientation-reversing automorphisms in this proposition are those
giving~(\ref{513:eqn}),~(\ref{514:eqn}), and~(\ref{515:eqn}).
\finedimo

\dimo{cosmetic:teo}
If $M=N$, then $\alpha,\beta\not\in\{\infty,-3,-2,-1,0\}$ because the
fillings listed in Theorem~\ref{one:cusp:teo} are all distinct. Then
$N(\alpha) = N(\beta)$ is hyperbolic. If $e(N(\alpha))=6$, by Corollary~\ref{e:Nalpha:cor}
we have $\{\alpha,\beta\} = \{-4,-3/2\}$, and the homeomorphism is
the orientation-reversing~(\ref{511right:eqn}).
Assume $e(N(\alpha))=5$, so $E(N(\alpha))=E(N(\beta))=\{\infty,-3,-2,-1,0\}$,
and let $f:N(\alpha)\to N(\beta)$ be the homeomorphism.
Note that $f$ acts on $\{\infty,-3,-2,-1,0\}$.
As in the previous proof, $\infty$ is the
only exceptional slope having distance $1$ from all other ones, so $f(\infty)=\infty$.
Then $f(n)$ is either $n$ or $-3-n$, so
either $N(\alpha,n)=N(\beta,n)$ for $n=-3,-2,-1,0$ or
$N(\alpha,n)=N(\beta,-3-n)$ for $n=-3,-2,-1,0$, but Table~\ref{two:cusps:table}
easily allows to conclude that this is impossible for $\alpha\neq\beta$.

Suppose now $M=N(\gamma)$. If $M(\alpha) = M(\beta)$ is not hyperbolic
then $\alpha,\beta$ are found by examining the repetitions in the
fillings in Table~\ref{two:cusps:table} and dismissing the pairs
related by the automorphisms of Proposition~\ref{repetitions:prop}.
The resulting pairs are precisely those of Table~\ref{cosmetic:table}.
Of these pairs, those giving as a filling the solid torus, which
is amphichiral, are both truly and
reflectively cosmetic. In the other cases the filling is chiral,
so the pair cannot be both truly and reflectively cosmetic.
For $\gamma\in\{-12/5,-6,-4/3\}$ we now explain how to
determine the type of the pair $\alpha,\beta$
by examining the fillings of $M(\alpha) = M(\beta)$.
In fact, using on $\partial M(\alpha)=\partial N(\gamma,\alpha)$ and
$\partial M(\beta)=\partial N(\gamma,\beta)$
the homology bases induced by $N$, we will have
$M\big(\alpha,\frac tu\big) = M\big(\beta,\frac {at+bu}{ct+du}\big)$
for some ${\tiny \matr abcd}\in\GL_2(\matZ)$,
and the pair $\alpha,\beta$ is truly or reflectively cosmetic
depending on whether ${\rm det} {\tiny \matr abcd}$ is $+1$ or $-1$.
For instance for $\gamma=-4/3$ we have from Table~\ref{closed:fill:table1}
\begin{eqnarray*}
M(-3,t/u)=N\big(-3,-4/3,t/u\big) & = & \seiftre {S^2}2152{t+u}u \\
M(-1,t/u)=N\big(-1,-4/3,t/u\big) & = & \seiftre {S^2}2152{t+3u}{-t-4u}
\end{eqnarray*}
whence ${\tiny \matr abcd}=\pm {\tiny \matr 47{-1}{-2}}$, so
$-3,-1$ are reflectively cosmetic.

To conclude we are left to show that $M(\alpha) = M(\beta)$
cannot be hyperbolic. Assume first $e(M(\alpha))>6$,
so $M(\alpha)$ appears in one of the Tables~\ref{one:cusp:table1}
to~\ref{one:cusp:table3}. By construction,
the manifolds listed in these tables are distinct and appear with all
their possible presentations as fillings of $N$. Whenever multiple presentations
exist one sees they are related by the symmetries of
Proposition~\ref{repetitions:prop}, so $\alpha$ and $\beta$ are equivalent.

Suppose now $e(M(\alpha))=6$, so $M(\alpha)$ appears in one of the
Tables~\ref{two:cusps:table1} to~\ref{one:cusp:table4}. Each
manifold in Table~\ref{one:cusp:table4} has a unique presentation
as a filling of $N$, so $M(\alpha)$ cannot appear there. If
$M(\alpha)$ appears in Table~\ref{two:cusps:table1}, so
$\gamma=1$, the only filling of $M(\alpha)$ admitting a
non-separating torus is $M(\alpha,1)$, and the topological type of
$M(\alpha,1)$ determines $\alpha$, whence $\alpha=\beta$, a
contradiction. If
$M(\alpha)$ appears in Table~\ref{two:cusps:table4}, the only
non-trivial symmetry of $E(M(\alpha))$ is $t/u\mapsto -t/u-4$, so
either $M(\alpha,t/u)=M(\beta,t/u)$ or
$M(\alpha,t/u)=M(\beta,-t/u-4)$ for all $t/u$. Looking at
Table~\ref{two:cusps:table4} itself, one easily sees that the
former equality is impossible for $\alpha\neq\beta$. If the latter
equality holds, using~(\ref{511left:eqn}) we get a slope $\beta'$
equivalent to $\beta$ such that
$M(\alpha,t/u)=M(\beta,-t/u-4)=M(\beta',t/u)$, giving
$\alpha=\beta'$ for what just said. Hence $\alpha$ and $\beta$ are
equivalent, a contradiction. A similar argument applies if
$M(\alpha)$ appears in Table~\ref{two:cusps:table2}
or~\ref{two:cusps:table3}.

If $e(M(\alpha))=5$ then $N(\gamma,\alpha,\delta) =
N(\gamma,\beta,h(\delta))$, where $h$ leaves $\{\infty,-3,-2,-1,0\}$
invariant, so either $h(\delta)=\delta$ or $h(\delta)=-3-\delta$.
With a little patience one sees from Table~\ref{closed:fill:table1}
that the latter possibility is actually
absurd, while the former implies that $\alpha=\beta$.
\finedimo

\subsection{Distance between exceptional slopes}
As stated in the Introduction, the maximal distance
$\Delta(X_0,X_1)$ of exceptional slopes giving manifolds of some
types $X_0$ and $X_1$
in $\{S,D,A,T,S^H,T^H,Z\}$
is known for most cases when both $X_0, X_1$
are not $Z$. It turns out that the partial
fillings of our chain-link complement $N$ allow to realize many
values of $\Delta(X_0,X_1)$, as summarized in
Table~\ref{maxima:table} and shown below.
For this reason we think that most of the lower bounds for $\Delta(X,Z)$
stated in the Introduction (and shown below)
should be optimal.

\begin{table}
\begin{center}
\begin{tabular}{c||c|c|c|c|c|c}
& $S$ & $D$ & $A$ & $T$ & $S^H$ & $T^H$ \\ \hline\hline
$S\phantom{\Big|}$ & {\bf 1} $(\geqslant 0)$ & {\bf 0} $(\geqslant -\infty)$ &
  {\bf 2} $(\geqslant -\infty)$ & {\bf 3} (3) &
  $\mathbf{-\infty}$ {\bf ?} ($\geqslant-\infty$) & {\bf 1} (1) \\ \hline
$D\phantom{\Big|}$ & & {\bf 1} (1) & {\bf 2} (2) & {\bf 2} (2) & --- & --- \\ \hline
$A\phantom{\Big|}$ & & & {\bf 5} (5) & {\bf 5} (5) & --- & --- \\ \hline
$T\phantom{\Big|}$ & & & & {\bf 8} (8) & {\bf 2} (2) & {\bf 3?} ($\geqslant 3$) \\ \hline
$S^H\phantom{\Big|}$ & & & & & {\bf 0} (0) & {\bf 1} (1) \\ \hline
$T^H\phantom{\Big|}$ & & & & & & {\bf 1} (1) \\
\end{tabular}
\end{center}
\nota{Known or conjectured values of $\Delta(X_0,X_1)$ shown in boldface,
and best lower bounds realized by fillings of $N$ shown between brackets.
Every known $\Delta(X_0, X_1)$ with $X_0,X_1\neq S$ is realized}
\label{maxima:table}
\end{table}

Table~\ref{families:table} 
describes some notable families of hyperbolic manifolds obtained
from $N$ by filling one or two cusps. It gives the topological type of the
fillings of these manifolds along the exceptional slopes
$\infty$, $-3$, $-2$, $-1$, and $0$. The table easily allows to prove most of
the lower bounds $\Delta(X_0, X_1)$
stated in Table~\ref{maxima:table} and in the Introduction. The remaining bounds are given
by the Whitehead link sister and figure-eight knot exteriors $N(1)$ and $N(1,2)$,
see Tables \ref{four:manifolds:table} and \ref{one:cusp:table1}.

\begin{table}
\begin{center}
\begin{tabular}{c||c|c|c|c|c}
\phantom{\Big|} & $\infty$ & $-3$ & $-2$ & $-1$ & $0$ \\ \hline \hline
\phantom{\Big|} $N(-2 + \frac 1n)$ & $D$ & $A, T$ & $D$ & $A$ & $A, T$ \\ \hline
\phantom{\Big|} $N(n, -4-n)$ & $T^H$ & $T$ & $Z$ & $Z$ & $S$ \\ \hline
\phantom{\Big|} $N\left(n, -4-n + \frac 1k\right)$ & $T^H$ & $T$ & $Z$ & $Z$ & $T^H$ \\ \hline
\phantom{\Big|} $N\left(-3 +\frac 1n, \frac{\pm(3n+1)+k(6n-1)}{\mp n-k(2n-1)}\right)$ &
  $T^H$ & $T$ & $Z$ & $S^H$ & $T$ \\ \hline
\phantom{\Big|} $N\left(-2+\frac 1n,\frac{kn\pm 1}{k-2(kn\pm1)}\right)$ & $S^H$ & $T$ & $T^H$ & $Z$ & $T$ \\ \hline
\phantom{\Big|} $N\left(-1 +\frac 1n, -1 -\frac 1n\right)$ & $S^H$ & $T^H$ & $Z$ & $Z$ & $T$ \\ \hline
\phantom{\Big|} $N\left(-\frac 52, \beta\right)$ & $T^H$ & $T$ & $T^H$ & $T^H$ & $T$

\end{tabular}
\end{center}
\nota{Exceptional fillings of some families $N(\alpha)$ or $N(\alpha, \beta)$
of cusped hyperbolic manifolds constructed from $N$.
The parameters $n, k\in\matZ$ and $\beta\in\matQ$ vary in the complement of a finite
set in $\matZ$ or $\matQ$.}
\label{families:table}
\end{table}

Some of the filled manifolds of Table~\ref{families:table} are noteworthy and
particularly easy to describe.
All boundary-reducible ($D$-type) manifolds in the first row are equal
to the solid torus. The reducible ($S$-type) manifold in the second row is
always $\matRP^3 \# \lens 31$.
The lens spaces ($T^H$-type) in row 2, row 3 last column, and row 6, are respectively
$\lens{n^2-4n-1}n$,
$\lens {6k-1}{2k-1}$, and
$\lens{4n^2+3}{2n^2+n+2}$.
Note that the only $S^H$-manifold is $S^3$.

Most of the families shown in Table~\ref{families:table} have already appeared elsewhere:
the first row consists of some complements of $1$-bridge braids in the solid torus, classified
by Berge~\cite{Be} and Gabai~\cite{Ga2},
the third row probably shows a family already found by Gordon~\cite{gordon}, and
the fourth row must consist of the Eudave-Mu\~noz knot complements~\cite{EuMu2} by a result of
Gordon and Luecke~\cite{GoLu2}. Finally, we do not know whether the fifth and sixth row give rise to
Berge knots, as conjectured in~\cite{Kirby}.

\end{document}